\documentclass[a4paper, 12pt]{article}
\usepackage{enumerate, theorem}
\usepackage{amsmath, amsfonts, amssymb}
\usepackage[height=22.5cm, width=15cm]{geometry}
\usepackage[all]{xy}
\usepackage{mathrsfs, yfonts}

\newtheorem{thm}{Th\'eor\`eme}[subsection]
\newtheorem{defn}[thm]{D\'efinition}
\newtheorem{cor}[thm]{Corollaire}
\newtheorem{prop}[thm]{Proposition}
\newtheorem{lemma}[thm]{Lemme}

\begin{document}

\title{\textbf{Sur certaines singularit\'es non isol\'ees d'hypersurfaces I.}}

\author{Daniel Barlet}
\date{6/12/05 .}
\maketitle
\markright{Sur certaines singularit\'es ... I}
\bigskip

\section*{Abstract.}
 
 \bigskip
 
 \noindent The aim of this fisrt part is to introduce, for an rather large class of hypersurface singularities with  1 dimensional locus, the analog of the Brieskorn lattice at the origin (the singular point of the singular locus). The main results are the finitness theorem for the corresponding (a,b)-module obtained via Kashiwara's constructibility theorem, and non torsion results for a plane curve singularity (not nessarily reduced) and for the suspension of such non torsion cases with an isolated  singularity.
 
 \bigskip
 
 \noindent {\bf Key words.} Hypersurface singularity, 1 dimensional singular locus, Brieskorn module,
 (a,b)-module, formal microlocal operators.
 
 \bigskip
 
 \noindent  {\bf Classification.} 32 S 05, 32 S 25, 32 S 40 .
 
\bigskip  

\section{Introduction}
\noindent L'objectif de cet article est de mettre en place, dans le cadre de fonctions  \`a lieu singulier de dimension 1, avec des hypoth\`eses assez restrictives mais donnant acc\`es \`a beaucoup d'exemples non triviaux, l'analogue de le th\'eorie de E. Brieskorn pour une fonction \`a singularit\'e isol\'ee (voir  [Br.70]).

\noindent Ce premier volet est centr\'e sur la construction de l'analogue du module de Brieskorn associ\'e \`a l'origine (qui est le point singulier du lieu singulier de \ $ f = 0 $) et les principaux r\'esultats obtenus sont
\begin{itemize}
\item{1)}  Le th\'eor\`eme de finitude \ref{finitude}.
\item{2)} Le th\'eor\`eme de non torsion pour les courbes planes non n\'ecessairement  r\'eduites \ref{P}.
\item{3)} La stabilit\'e de l'absence de torsion par suspension avec une fonction \`a  singularit\'e isol\'ee  (Proposition \ref{T.S.2}.)
\end{itemize}

\smallskip

\noindent Le th\'eor\`eme de finitude nous permet d'attacher \`a une fonction v\'erifiant notre hypoth\`ese \ $(HI)$ \ un (a,b)-module \`a l'origine. Malheureusement, le calcul du rang de ce (a,b)-module, qui donne la dimension du $n-$i\`eme groupe de cohomologie de la fibre de Milnor de \ $f$ \`a l'origine, fait intervenir la torsion \'eventuellement pr\'esente. Quand cette torsion est nulle, on obtient une jolie formule g\'en\'eralisant celle de\\ J. Milnor pour le cas d'une singularit\'e isol\'ee.

\noindent Nous donnons  alors des conditions simples pour assurer l'annulation de cette torsion. D'o\`u l'int\'er\^et des r\'esultats 2) et 3)  de non torsion.

\noindent Le second volet (voir [B.04 b]) s'attaquera \`a une version filtr\'ee du ph\'enom\`ene d'interaction de strates cons\'ecutives \'etudi\'e dans [B.91]. Ceci n\'ecessitera au pr\'ealable la g\'en\'eralisation au cas de la valeur propre 1 les r\'esultats de loc. cit.

\noindent Cette g\'en\'eralisation est d\'ecrite  dans [B.04 a], mais n\'ecessite l'\'etude d'une situation \`a trois strates ; on consultera  l'introduction de [B.05] pour des commentaires plus d\'etaill\'es sur les interactions de strates.

\section{ Pr\'e-(a,b)-Modules.}

\subsection{G\'en\'eralit\'es.}

\begin{defn}\label{pre}

\noindent Soit  E  un espace vectoriel complexe muni d'endomorphismes \ $a$ \ et  \ $b$. On pose 
$$ B(E) =  \bigcup_{m\in \mathbb{N}}Ker\,  b^m  \quad {\rm et} \quad A(E) = \lbrace x \in E \  / \  \mathbb{C}[b].x \subset  \bigcup_{m\in \mathbb{N}} Ker\, a^m \ \rbrace $$

\noindent On dira que \ $E$ \ est un {\bf pr\'e-(a,b)-module} lorsque  les conditions suivantes sont v\'erifi\'ees 
\begin{enumerate}[i)]
\item  \quad $ a.b - b.a = b^2  \ . $

\item  Pour tout \ $ \lambda \in \mathbb{C}^*, \ b - \lambda $ \ est bijectif  dans \ $E$.

\item   \ $ \exists N \in \mathbb{N} \ | \quad  a^N .A(E) = 0 $.
\item  \  $ B(E) \subset A(E) \ . $
\item   \ $  \bigcap_{m\in \mathbb{N}} b^m (E) \subset A(E) \ . $
\item   Le noyau et le conoyau de \ $b$ \ sont de dimensions finies sur \ $\mathbb{C} . \hfill \Box $
\end{enumerate}

\noindent On dira que \ $E$ \ est {\bf sans torsion} (en fait sans \ $b-$torsion) si on a de plus  \ $ Ker\, b = 0 $ \ ce qui \'equivaut \`a \ $ B(E) = 0  .$  Dans ce cas la condition \ $ B(E) \subset A(E) $ \  devient triviale. Mais il resultera du lemme \ref{A=B} que dans ce cas on a \'egalement \ $ A(E) = 0 $. 
\end{defn}

\bigskip

\noindent  {\bf Remarques. }
\begin{itemize}
\item{1)} Un (a,b)-module, c'est \`a dire un \ $\mathbb{C}[[b]]-$module libre de type fini muni d'un endomorphisme \ $\mathbb{C}-$lin\'eaire \ $a$ \ v\'erifiant la relation de commutation  i),
est un pr\'e-(a,b)-module. En effet l'injectivit\'e de \ $b$ \ donne la nullit\'e de \ $B(E)$; donc  la condition iv) est triviale. La \ $b-$compl\'etion, qui implique la \ $b-$s\'eparation, donne  la condition v) et la condition ii).
La condition  vi)  est \'egalement \'evidente. Il reste \`a montrer la  condition iii). Montrons directement que l'on a \ $A(E) = 0 $ \ dans ce cas.

\noindent Soit \ $ x \in A(E)$ \ v\'erifiant \ $ a.x = 0 $. Alors on doit avoir \ $ a^n.bx = 0 $ \ pour \ $ n \in \mathbb{N} $ \ assez grand.  Mais si on a \ $ a.x = 0 $ \ et \ $ a^n.bx = 0 $ \ on en d\'eduit que \ $ a^{n-1}.(ba + b^2).x = 0 $ \ d'o\`u \ $ a^{n-1}.b^2 x = 0$ \ puis \ $ a^{n-p}.b^p x = 0 $ \ pour \ $ p = 1, \dots , n$. On a donc \ $ b^n.x = 0 $ \ d'o\`u \ $ x = 0 $ \ puisque \ $b$ \ est injective. Comme, par d\'efinition, \ $a$ \ est nilpotente sur \ $A(E)$ \ on en conclut que \ $ A(E) = 0 $ \ pour un (a,b)-module. La condition iii) est donc trivialement v\'erifi\'ee dans ce cas.
\item{2)} Gr\^ace \`a la condition ii), un pr\'e-(a,b)-module est un module sur le localis\'e  \ $ \mathbb{C}[b]_0 $ \  de l'anneau \ $ \mathbb{C}[b] $ \  par rapport \`a l'id\'eal maximal engendr\'e par \ $b$.
\item{3)} On remarquera que, par d\'efinition,  $A(E)$ \ est stable par \ $b$ \ et que\ $x\in A(E)$ \ si  et seulement  si\ $ \forall p\in \mathbb{N} , \exists n\in \mathbb{N}\  /\  a^n.b^p.x = 0 $. On en d\'eduit que \ $A(E)$ \ est stable par \ $a$ \ car si on a \ $ a^N.b^p.x = a^N.b^{p+1}.x = 0 $ \ alors \\
 $$ a^N.b^p.ax = a^N.(a.b^p - p.b^{p+1}).x = 0 . $$
\end{itemize}

\begin{lemma} \label{A=B}

\noindent Si \ $E$ \ est un pr\'e-(a,b)-module, alors \ $A(E)$ \ est un \ $\mathbb{C}-$espace vectoriel de dimension finie. De plus on a l'\'egalit\'e \ $ B(E) = A(E) $.
\end{lemma}

\bigskip

\noindent \textit{\underline{Preuve}.} Comme \ $A(E)$ \ est stable par \ $a$ \ et \ $b$ \ l'\'egalit\'e \ $  a^N .A(E) = 0 $ \ donne \ $ b^{2N}.A(E) = 0 $ . En effet la relation de commutation \  i) \  implique la formule
$$ N!\, b^{2N} = \sum_{j=0}^{N} (-1)^j \big(^j_N\big).b^j a^N b^{N-j}  \ . $$
Elle est \'etablie dans [B.95] p.24. 

\noindent  L'\'egalit\'e \ $  B(E) = A(E) $ \ s'en d\'eduit en constatant que dans le quotient \ $ A(E) \big/ B(E) $ \  l'endomorphisme \ $b$ \ est \`a la fois injectif, puisque \ $ bx \in B(E) $ \ implique \ $ x \in B(E)$, et nilpotent, puisque \ $ b^{2N}.A(E) = 0 $.

\noindent Montrons que l'espace vectoriel  \ $  B(E) = A(E) $ \ est de dimension finie. On a pour chaque \ $ \nu \in \mathbb{N} $ \ la suite exacte de \ $\mathbb{C}-$espaces vectoriels  :
$$ 0 \rightarrow Ker\, b \longrightarrow Ker\, b^{\nu + 1} \overset{b}{\longrightarrow} Ker\, b^{\nu} \qquad  ; $$
elle donne, par r\'ecurrence, grace \`a la finitude de \ $ Ker\, b $, la finitude de la dimension de \ $ Ker\, b^{\nu}  \quad  \forall \nu \in \mathbb{N} $. Mais on a \ $ A(E) \subset  Ker\, b^{2N} $.\ $\hfill \blacksquare $

\bigskip

\newpage

\noindent {\bf Remarques. }
\begin{enumerate}[1)]
\item La preuve des assertions \ $ b^{2N}A(E) = 0 $ \ et \ $B(E) = A(E) $ \ dans le lemme \ref{A=B} n'utilise pas la condition de finitude vi) pour  \ $E$.

\item Le lemme pr\'ec\'edent donne \'egalement, sans utiliser la condition vi) pour \ $E$, la \ $b-$s\'eparation, c'est \`a dire la nullit\'e de \ $\cap_{m\geq 0} b^m(E) $. En effet, si \ $x \in \cap_{m\geq 0} b^m(E) $ \ on peut trouver \ $y \in E $ \ v\'erifiant \ $b^{2N}.y = x $. Mais comme \ $ x \in B(E) $ \ on a \'egalement \ $ y \in B(E) = A(E) $, et donc \ $ b^{2N}.y = x = 0.$ Les propri\'et\'es i) \`a v) seules, suffisent donc \`a assurer la \ $b-$s\'eparation de \ $E$ , et donc son injection dans son compl\'et\'e \ $b-$adique, sans avoir \`a quotienter par \ $A(E) = B(E)$ , c'est \`a dire par la \ $b-$torsion.
Ceci s'applique, en particulier, \`a la situation de la proposition \ref{fonction}, ce qui couvre tous les cas issus de singularit\'es d'une fonction holomorphe (arbitraire) que nous consid\`ererons.

\item  La \ $a-$torsion de \ $E$, qui est le sous-espace vectoriel \ $\tilde{A}(E) : = \bigcup_{m\in \mathbb{N}} Ker\, a^m $, peut \^etre strictement plus gros que \ $ A(E)$ \ et ne pas \^etre stable par \ $b$ \ m\^eme dans le cas d'un (a,b)-module. Prendre par exemple le (a,b)-module \ $ \mathbb{C}[[b]]-$libre de rang 1  et de g\'en\'erateur not\'e  \ $e$ \ dont l'application \ $a$ \ est d\'efinie par \ $ a.e = 0 $ \ (alors \ $ab^n .e = n.b^{n+1}.e \quad  \forall n\in\mathbb{N}$). 

\item \`A titre d'exercice, le lecteur pourra montrer que pour un pr\'e-(a,b)-module r\'egulier (voir \ref{reg}) l'espace vectoriel \ $\tilde{A}(E)$ \ est toujours de dimension finie (il suffit en fait de traiter le cas d'un (a,b)-module \`a p\^ole simple).

\item Nous montrerons au lemme \ref{positiv.} que pour les pr\'e-(a,b)-modules associ\'es aux syst\`emes de Gauss-Manin d'un germe de  fonction holomorphe on a toujours l'\'egalit\'e \ $ \tilde{A}(E) = A(E) $ \ grace au th\'eor\`eme de positivit\'e de B. Malgrange [M.74].

\end{enumerate}

\begin{prop}\label{(a, b)}

\noindent Soit \ $E$ \ un pr\'e-(a,b)-module. Le compl\'et\'e \ $b-$adique \ $\mathcal{L}(E)$ \ du quotient  \ $E \big/ A(E) $ \ est un (a,b)-module (voir plus haut ou bien [B.93], [B.95] ou [Be. 01]). L'application (a,b)-lin\'eaire naturelle 
$$ E \longrightarrow \mathcal{L}(E) $$
a pour noyau \ $B(E) = A(E) $ ; elle est continue  et d'image dense pour la topologie \ $b-$adique.
\end{prop}

\bigskip

\noindent \textit{\underline{Preuve}.} D'apr\`es ce qui pr\'ec\`ede, \ $b$ \ est injective sur \ $ E\big/ A(E) $ \ qui est s\'epar\'e pour la filtration \ $b-$adique d'apr\`es la condition  v). Pour conclure il suffit de voir que le conoyau de \ $b$ \ agissant sur \ $ E\big/ A(E) $ \ est de dimension finie sur \ $\mathbb{C} $ . Ceci est clair puisque c'est un quotient de \ $ E/bE $ .  \ $\hfill \blacksquare $

\begin{defn}
Nous appelerons {\bf rang} du pr\'e-(a,b)=module \ $E$,  not\'e \ $rg(E)$, le rang (comme \ $\mathbb{C}[[b]]-$module libre ) du (a,b)-module \ $\mathcal{L}(E)$ \ qui lui est associ\'e.
\end{defn}

\noindent On remarquera que si \ $\delta : = \dim Ker\, b $, on a \ $ rg(E) = \dim E/b.E - \delta $, puisque \ $B(E)$ \ est de dimension finie.

\subsection{R\'egularit\'e et produit tensoriel.}

\begin{defn}\label{reg}

\noindent Nous dirons qu'un pr\'e-(a,b)-module \ $E$ \ est {\bf local}, respectivement {\bf \`a p\^ole simple}, respectivement  {\bf r\'egulier}, quand il v\'erifie :
\begin{itemize}
\item{\bf Local } : \ $ \exists l \in \mathbb{N}^* \quad {\rm tel \ que} \quad   a^l.E \subset b.E \ .$
\item{\bf \`A p\^ole simple} :   \ $ l = 1 $ \ dans la condition pr\'ec\'edente, c'est \`a dire  \ $ a.E \subset b.E \ . $ 
\item{\bf R\'egulier} : \ $ \exists k \in \mathbb{N}^*   \quad {\rm tel \ que} \quad  a^k.E \subset \underset{j\in [0,k-1]}{\sum} b^{k-j}.a^j .E \ .  \hfill \Box $
\end{itemize}
\end{defn}

\bigskip

\noindent {\bf Remarques.}
\begin{itemize}
\item  On a, bien sur, "p\^ole simple" \ $\Rightarrow$ \ "r\'egulier"  \ $\Rightarrow$ \ "local". 
\item  Il est \'equivalent de demander que \ $E$ \ soit local (resp. r\'egulier), ou bien que \ $ E/B(E) $ \ soit local  (resp. r\'egulier), ou encore que \ $\mathcal{L}(E)$ \ le soit. 

\end{itemize}

\begin{defn}\label{prod.tens. 0}
\noindent Soient \ $E$ \ et \ $F$ \ deux pr\'e-(a,b)-modules ; alors le produit tensoriel
$$ E \otimes_{\mathbb{C}[b]_0} F $$
 muni de l'application \ $\mathbb{C}-$lin\'eaire  \ $ a : =  (a_E \otimes 1_F + 1_E \otimes a_F ) \ $
sera appel\'e "produit tensoriel" de \ $E$ \ et \ $F$ \ et not\'e simplement \ $ E \otimes_{a,b} F$ \ ou plus simplement \ $ E \otimes F $ \ quand il n'y a pas d'ambiguit\'e.  $\hfill \Box $
\end{defn}

\begin{prop}\label{prod.tens. 1}

\noindent Soient  \ $E$ \ et \ $F$ \ deux pr\'e-(a,b)-modules locaux ;  alors  \ $ E \otimes_{a,b} F $ \  est \'egalement un pr\'e-(a,b)-modules local. Si de plus \ $E$ \ et \ $F$ \ sont r\'eguliers (resp. \`a p\^oles simples) alors \ $ E \otimes_{a,b} F $ \ est r\'egulier (resp. \`a p\^ole simple).
\end{prop}

\bigskip

\noindent \textit{Preuve.} La condition  i) se v\'erifie facilement.
\noindent  La condition  ii) est  imm\'ediate. Pour montrer la  condition iv), commen{\c c}ons par prouver que l'on a 
$$ B(E \otimes F) = B(E)\otimes F + E\otimes B(F) \ . $$
L'inclusion \ $B(E)\otimes F + E\otimes B(F) \subset  B(E \otimes F)$ \ est claire. Montrons l'inclusion oppos\'ee. Soit \ $ z \in  B(E \otimes F) $. Posons \ $ z  = \sum_{i\in I} x_i \otimes y_i $ \ o\`u \ $ I$ \ est fini, et notons par \ $E_1$ \ et \ $F_1 $ \ les sous-$\mathbb{C}[b]_0-$modules engendr\'es par les \ $(x_i)_{i\in I}$ \ et les \ $(y_i)_{i\in I}$ \ respectivement. Comme, par d\'efinition, \ $E_1$ \ et \ $F_1$ \ sont des \ $\mathbb{C}[b]_0-$modules de type finis, ils sont somme directe de leur torsion et d'un module libre de type fini. On en conclut ais\'ement que l'on a \ $ B(E_1 \otimes F_1) \simeq \big(B(E_1) \otimes F_1\big) +  \big( E_1 \otimes B(F_1)\big) $ \ ce qui nous donne \ $ z\in B(E)\otimes F + E\otimes B(F) . $

\noindent On en conclut que \ $ B(E \otimes F) \subset A(E \otimes F)$ \ en remarquant que si \ $ N \in \mathbb{N} $ \ est assez grand pour que l'on ait \  $ a^N.B(E) = 0 = b^N.B(E) $ \ ainsi que  \  $ a^N.B(F) = 0 = b^N.B(F) $ \ et si \ $M \in \mathbb{N} $ \ est assez grand pour que l'on ait \ $ a^M.E \subset b.E $ \ ainsi que \ $a^M.F \subset b.F $ \  on aura
 $$a^{M.N + N}.\big(B(E)\otimes F + E\otimes B(F)\big) = 0  $$
puisque
$$ a^{M.N + N}.\big(B(E)\otimes F\big) \subset \sum_{p\in [0,N]} a^p.B(E)\, \otimes\, a^{M.N}.F  \subset B(E)\,\otimes \, b^N.F = b^N.B(E)\, \otimes\,  F  .$$
Comme \ $B(E \otimes F)$ \ est stable par \ $b$ \ l'inclusion  \ $ B(E \otimes F) \subset A(E \otimes F)$ \   est v\'erifi\'ee.

\noindent Maintenant, grace \`a ce qui pr\'ec\`ede, pour \ $ N \in \mathbb{N} $ \  assez grand,
\ $ B(E \otimes F)$ \ sera un quotient de 
$$ \big( B(E)\, \otimes \, (F/b^N.F)\big) \oplus \big(E/b^N.E)\, \otimes \, B(F) \big)   $$
ce qui prouve la finitude de \ $ Ker\, b \subset B(E \otimes F)$.

\noindent La finitude sur \ $\mathbb{C}$ \ de \ $ Coker\, b $ \ dans  \ $ E \otimes_{a,b} F $ \ est \'el\'ementaire.

\noindent Montrons maintenant que  \ $E \otimes_{a,b} F $ \  v\'erifie la condition iii)  de la d\'efinition \ref{pre}. D'apr\`es ce qui pr\'ec\`ede  on a un isomorphisme (a,b)-lin\'eaire
$$ E \otimes_{a,b} F \big/ B( E \otimes_{a,b} F ) \simeq E/B(E) \otimes_{a,b} F/B(F)  \ .$$
Consid\'erons maintenant \ $ z \in A(E\otimes F) $ \ et supposons que \ $a.z \in B(E \otimes F) $. Comme \ $ b.z \in \bigcup_{m\in \mathbb{N}} Ker\, a^m $ \ il existe \ $m \in \mathbb{N} $ \ tel que 
$ a^m.b.z = 0$. On obtient alors \ $a^{m-i}b^{i+1}.z \in  B(E \otimes F), \forall  i\in [0,m] $, d'o\`u  \ $ b^{m+1}.z \in B(E \otimes F)$ \ et donc \ $z \in  B(E \otimes F)$.

\noindent Ceci montre que sur le quotient \ $ A( E \otimes_{a,b} F ) \big/ B( E \otimes_{a,b} F ) $ \ l'endomorphisme nilpotent \ $a$ \ est injectif. Donc \ $ A( E \otimes_{a,b} F ) = B( E \otimes_{a,b} F )$, et \ $a^{M.N + N}.A(E\otimes F) = 0$.

\noindent Il nous reste seulement \`a prouver la condition v) de la d\'efinition \ref{pre}. Mais \ $ E/B(E)$ \ et \ $F/B(F) $ \ sont \ $b-$s\'epar\'es  ainsi que 
$$E\otimes F/B(E\otimes F) \simeq E/B(E) \otimes F/B(F) \subset \mathcal{L}(E) \otimes \mathcal{L}(F) .$$
On en d\'eduit la condition v) pour \ $ E \otimes_{a,b} F$.

\noindent Le cas de pr\'e-(a,b)-modules  \`a p\^ole simple est \'evident. 

\noindent Pour traiter le cas de pr\'e-(a,b)-modules r\'eguliers, nous pouvons supposer  que ce sont des (a,b)-modules, d'apr\`es la remarque qui suit la d\'efinition \ref{reg}. Soit \ $K$ \ le corps des fractions de l'anneau \ $\mathbb{C}[[b]] $ . La condition de r\'egularit\'e pour un (a,b)-module  \ $G$ \  \'equivaut \`a la finitude dans 
$$ G \otimes_{\mathbb{C}[[b]] } K $$
du sous-\ $\mathbb{C}[[b]]-$module \ $ \sum_{m\geq0} (b^{-1}.a)^m .G $ \  (voir [B.93] p. 18 ).

\noindent Mais on a pour chaque \ $m \in \mathbb{N} $ 
$$ (b^{-1}.a)^m \big(E \otimes F\big) \subset \sum_{j\in [0,m]} \big((b^{-1}.a)^j E \otimes (b^{-1}.a)^{m-j} F \big)  \ . $$
On en d\'eduit la r\'egularit\'e de \ $ E \otimes_{a,b} F $ . $\hfill \blacksquare $

\subsection{Cas des pr\'e-(a,b)-modules associ\'es \`a des singularit\'es de fonctions holomorphes.}

\noindent Nous allons montrer que la cohomologie du complexe de De Rham 
\ $ ((Ker\,df)^{\bullet}, d)$ \  restreint \`a \ $f^{-1}(0)$, associ\'ee \`a une fonction holomorphe non constante, v\'erifie toujours les propri\'et\'es i) \`a v) de la d\'efinition d'un pr\'e-(a,b)-module. 

\begin{prop}\label{fonction}
\noindent Soit \ $ f : (\mathbb{C}^{n+1},0) \rightarrow (\mathbb{C},0) $ \ un germe non nul de fonction holomorphe. Consid\'erons pour \ $ p \in [1,n+1] $ \ le \ $\mathbb{C}-$espace vectoriel 
 $$ \mathcal{H}^p_0 : = \big((Ker\, df)^p \cap Ker\, d \big) \big/ d((Ker\, df)^{p-1})_0 $$
muni des \ $\mathbb{C}-$endomorphismes \ $a$ \ et \ $b$ \ d\'efinis respectivement par la multiplication par \ $f$ \ et \ $ df\wedge d^{-1}$.

\noindent Alors \ $ \mathcal{H}^p_0$ \ v\'erifie les conditions i) \`a v) de la d\'efinition \ref{pre}.
\end{prop}

\noindent \textit{\underline{Preuve}. }Soit \ $ du\in (Ker\, df)^p $ \ ; on a \ $ a[du] = [f.du] $ \ et \ $ b[du] = [df\wedge u] $. Donc
 $$ b(a + b)[du] = b([d(f.u)]) = [df \wedge f.u] = a([df\wedge u]) = a(b([du]))  \ , $$
ce qui prouve notre premi\`ere assertion.

\noindent Comme on a 
$$ d\big(\lambda .e^{-\frac{1}{\lambda}.f }.u \big) = -e^{-\frac{1}{\lambda}.f }\big(df\wedge u - \lambda du \big)  $$ 
l'annulation de \ $ b[du] - \lambda.[du] $ \ qui se traduit par l'existence de \ $ v \in (Ker\, df)^{p-1} $ \ v\'erifiant 
$$ e^{\frac{1}{\lambda}.f }. d\big(\lambda .e^{-\frac{1}{\lambda}.f }.u \big) = dv $$
ce qui donne 
$$  d\big(e^{-\frac{1}{\lambda}.f }.(\lambda.u - v) \big) = 0  .$$

 Le lemme de De Rham holomorphe assure alors l'existence de \ $ \xi \in \Omega^{p-2}_0 $ \ tel que l'on ait
$$ d\xi = e^{-\frac{1}{\lambda}.f }.(\lambda.u - v)   \ .$$
On en d\'eduit que
$$  d\big( e^{\frac{1}{\lambda}.f }.\xi\big) = \lambda.u - v  + \frac{1}{\lambda} e^{\frac{1}{\lambda}.f }.df\wedge\xi $$
ce qui donne la nullit\'e de \ $ \lambda.[du] $ \ et donc de \ $[du]$ \ puisque \ $\lambda \not= 0 .$

\noindent Pour voir la surjectivit\'e de \ $ b - \lambda $ \ il suffit de r\'esoudre l'\'equation 
$$ df \wedge \xi - \lambda .d\xi = du  $$
o\`u \ $ du \in (Ker\, df)^p $ \ est donn\'e. Comme la p-forme holomorphe \ $  -\frac{1}{\lambda}.e^{\frac{-1}{\lambda}.f }.du $ \ est \ $d-$ferm\'ee, le lemme de De Rham holomorphe assure alors l'existence de \ $ \eta \in \Omega^{p-1}_0 $ \ v\'erifiant
$$ d\eta = -\frac{1}{\lambda}.e^{\frac{-1}{\lambda}.f }.du  \ .$$
On v\'erifie alors imm\'ediatement que \ $ \xi : = e^{\frac{1}{\lambda}.f}.\eta $ \ est solution de notre \'equation\footnote{Pour \ $p=1$ \ on choisit \ $u(0) = 0$. Comme \ $(Ker\,df)^0= 0 $, on aura ensuite \ $v = 0$ \ et \ $ u = C.\lambda^{-1}.e^{\frac{1}{\lambda}.f }$. Alors \ $u(0) = 0$ \ donne \ $C = 0$ \ et donc \ $ u \equiv 0$.}.   

\noindent La condition iii) est  cons\'equence du r\'esultat classique d'A. Grothendieck  [G.66]\footnote{On peut obtenir directement ce point via le cas \`a croisements normaux en utilisant le th\'eor\`eme de d\'esingularisation d'Hironaka [H.64].} . La condition  iv)  est cons\'equence du fait que la \ $b-$torsion de \ $\mathcal{H}^p_0$ \ est contenue dans la \ $a-$torsion. Ceci r\'esulte du fait que, quitte \`a localiser en \ $a$ \ (c'est \`a dire \`a consid\'erer des formes m\'eromorphes \`a p\^oles dans \ $f=0$ \ ) on a \ $ Ker\, df^{\bullet} = Im\, df^{\bullet} $ . 
\begin{align*}
{\rm Alors} \quad \quad df \wedge \xi  & =  & d\eta \quad {\rm avec} \quad \quad & df\wedge \eta & = &\  0 \quad {\rm donnera}\ & \eta & =\  df \wedge \zeta \ ;\\
{\rm donc} \quad df\wedge (\xi + d\zeta) &  = & 0 \quad  {\rm et\  donc}\quad & \xi + d\zeta &  =  & \ df \wedge \gamma \ .\  {\rm Alors} \ & d\xi & =  -df\wedge d\gamma \ .
\end{align*}
La condition  v)  est une cons\'equence imm\'ediate du th\'eor\`eme de positivit\'e de B. Malgrange [M.74]. $\hfill \blacksquare $

\begin{lemma}\label{positiv.}
\noindent Pla{\c c}ons-nous dans la situation de la proposition  \ref{fonction}  et supposons la fonction \ $f$ \ r\'eduite. Pour tout \ $ p\in \mathbb{N}\ ,  p \geq 2\  ,  E : =\mathcal{H}^p_0 $ \ v\'erifie \ $ \tilde{A}(E) = A(E) $.
\end{lemma}

\bigskip

\noindent \textit{\underline{Preuve}.} Il s'agit simplement de voir que \ $ \tilde{A}(E) $ \ est stable par \ $b$. Ceci r\'esulte imm\'ediatement du th\'eor\`eme de positivit\'e de B. Malgrange\footnote{L'assertion \ "$f$ \ r\'eduite" correspond \`a la preuve pr\'ecis\'ee dans l'Appendice de [B.84] (note au bas de la page 106 ).} qui affirme que pour une forme holomorphe \ $\omega \in (Ker\, df)^p$ \ et pour \ $\gamma_s \subset f^{-1}(s)$\  une famille horizontale multiforme de \ $(p-1)-$cycles compacts,  la fonction holomorphe multiforme
$$ s \rightarrow \int_{\gamma_s} \frac{\omega}{df} $$
 tend vers \ $0$ \ quand \ $s$ \ tend vers \ $0$ \ le long de tout rayon issu de l'origine (ici on utilise \ $ p-1 \geq 1).$

\noindent Alors \ $\omega$ \ induit une classe de \ $a-$torsion dans \ $E$ \ si et seulement si pour tout choix de la famille \ $\gamma_s $ \ on trouve une int\'egrale nulle. Comme l'op\'eration \ $b$ \ correspond \`a prendre une primitive (qui doit donc \^etre "nulle en \ $0$"), l'assertion \ $ \tilde{A}(E) = A(E) $ \ s'en d\'eduit. $ \hfill \blacksquare $ 

\begin{cor}

\noindent Si \ $f$ \ est \`a singularit\'e isol\'ee \`a l'origine, \ $\mathcal{H}^{n+1} $ \ est un pr\'e-(a,b)-module sans torsion, dans lequel \ $a$ \ est injective.
\end{cor}

\bigskip

\noindent  Ceci donne une nouvelle preuve du  r\'esultat de Sebastiani [S.70].

\bigskip

\noindent \textit{\underline{Preuve}.} V\'erifions d\'eja la condition vi)  de la d\'efinition \ref{pre}. La dimension finie de \ $ Coker\, b $ \ est imm\'ediate vue que ce noyau s'identifie au quotient \ $ \Omega^{n+1}_0\big/ df \wedge \Omega^n_0 $ . Montrons que \ $Ker\, b = 0 $. Soit \ $ \omega = d\xi $ \ o\`u \ $ \xi \in \Omega^n_0 $ \ v\'erifiant \ $ df\wedge \xi = df \wedge d \eta $ \ (rappelons que l'on a \ $ Ker\, df^n = df \wedge \Omega^{n-1} $ \ puisque \ $f$ \ est \`a singularit\'e isol\'ee). Alors \ $ \xi = d\eta + df \wedge \zeta $ \ et on a donc \ $ [\omega] = 0 $ \ dans \ $ \mathcal{H}^{n+1} .$

\noindent Il suffit alors d'appliquer le lemme pr\'ec\'edent pour conclure, puisque l'on a vu que \ $0 =  B(E) = A(E). \hfill \blacksquare $

\section{Le th\'eor\`eme de finitude.}

\subsection{L'hypoth\`ese \ $ (HI)$.}

\noindent Soit $\;\tilde{f}:(\mathbb{C}^{n+1},0) \longrightarrow (\mathbb{C},0)\;$ un germe non constant
de fonction holomorphe et soit $\;f:X \longrightarrow D\;$ un repr{\'e}sentant de Milnor
de $\;\tilde{f}.\;$ Nous ferons les hypoth{\`e}ses suivantes:
\smallskip

\begin{itemize}
\item [{\it HI~a)}] Le lieu singulier $\;S:=\{x \in X\big/df_x=0\}\;$ est une courbe
            contenue dans $\;Y:=f^{-1}(0),\;$ dont chaque composante irr{\'e}ductible contient
            l'origine et est non singuli{\`e}re en dehors de $\;0$.
\item [{\it HI~b)}] En chaque point  $\;x\;$ de $\;S-\{0\}\;$ il existe un germe
            en $\;x\;$ de champ de vecteur holomorphe $\;V_x\;$, non nul en $\;x,\;$ tel que
            $\;V_x \cdot f \equiv 0$.
\end{itemize}
\smallskip

\noindent L'hypoth{\`e}se {\it HI~b)} est assez restrictive puisqu'elle implique que le long
de \\ $\;S-\{0\}\;$ la singularit{\'e} $\;\{f=0\}\;$ est une d{\'e}formation localement triviale
de la singularit{\'e} hyperplane transverse (qui est une singularit{\'e} isol{\'e}e de $\;
\mathbb{C}^{n})$.

\noindent Cependant il est facile de voir que cette hypoth{\`e}se est toujours v{\'e}rifi{\'e}e pour
$\;n=1\;$ (courbes planes r{\'e}duites ou non ) et qu'il y a beaucoup d'exemples en
dimensions sup{\'e}rieures ( voir le paragraphe 4.2.)

\bigskip

\noindent Le r\'esultat fondamental de ce  paragraphe est le

\begin{thm}\label{finitude}
 Sous l'hypoth\`ese (HI), l'espace vectoriel  
$$ \mathcal{H}^{n+1}: = \Omega^{n+1}_0 \big/ d(Ker\, df)^n_0  $$ 
muni des endomorphismes \ $a$ \ et \ $b$ \ donn\'es respectivement par multiplication par \ $f$ \ et par \ $ df\wedge d^{-1}$ \  est un pr\'e-(a,b)-module r\'egulier.
\end{thm}

\noindent Compte tenu des r\'esultats du paragraphe 2  il s'agit  maintenant de prouver les conditions de finitudes vi)  de la d\'efinition \ref{pre}. Ceci utilisera, entre autres, le th\'eor\`eme de constructibilit\'e de M. Kashiwara [K.75].  

\subsection{ L'id{\'e}al $\;\widehat{J(f)}$.}

\noindent Dans la situation pr{\'e}cis{\'e}e ci-dessus, introduisons l'id{\'e}al $\;\widehat{J(f)}\;$
de $\;{\cal O}_X\;$ form{\'e} des germes de fonctions holomorphes dont la restriction
{\`a} $\;X-\{0\}\;$ est dans $\;J(f) $ \ l'id\'eal jacobien de \ $f$. Si $\;i:X-\{0\} \hookrightarrow X\;$ est 
l'inclusion, on a, d{\`e}s que $\;n \ge 1,\;$ d'apr{\`e}s Hartogs
$$\widehat{J(f)} \simeq i_* i^* \big(J(f)\big).$$

\begin{defn}
Le faisceau $\;\widehat{J(f)}\big/J(f)\;$ est coh{\'e}rent et concentr{\'e} en $\;0$. En effet, le choix d'un {\'e}l{\'e}ment de volume sur $\;X\;$ donne un isomorphisme de faisceaux 
$$\widehat{J}(f)\big/J(f) \overset{\sim}{\longrightarrow} \underline{H}^0_{\{0\}} \big(
\Omega^{n+1}_X \big/ df_{\wedge} \Omega^n_X\big).$$

\noindent La dimension finie de son germe en $\;0\;$ sera not{\'e}e \ $\mu(f)$.
\end{defn}

\bigskip

\noindent {\bf Remarque.}

\noindent Pour $\;f\;$ a singularit{\'e} isol{\'e}e en $\;0,\;$ on a $\;\widehat{J(f)}={\cal O}_X\;$
et donc  \ $\mu = {\rm dim}_{\mathbb{C}} \; {\cal O}\big/J(f)$ \  ce qui est bien conforme
{\`a} la d{\'e}finition du nombre de Milnor dans ce cas.

\subsection{Le $\; {\cal D}-$module $\; {\cal M} = {\cal D} \big/ {\cal J}$.}

\noindent Soit Ann$(f)$ le sous-faisceau du faisceau $\;T_X\;$ des champs de vecteurs
holomorphes sur $\;X$, form\'e des germes de champs de vecteurs qui annulent  \ $f$. Ce sous-faisceau est coh{\'e}rent car un choix de coordonn{\'e}es sur $\;X\;$ montre qu'il est isomorphe via
$$\displaystyle  \sum^n_{i=0} a_{i} \; \frac{\partial}{\partial x_{i}} \longrightarrow \sum^n_{i=0}
a_{i} \; \widehat{\overset{i}{dx}}$$
 au noyau  \ $ (Ker\, df)^n$ \ du morphisme $\;{\cal O}_X -$lin{\'e}aire 
$$ df^n\wedge : \Omega^n_X \longrightarrow \Omega^{n+1}_X \quad \hbox{d{\'e}fini par} \quad
\alpha \longrightarrow df_{\wedge} \alpha.$$

\noindent La condition $\;(HI)\;$ implique que, si $\;V_1 \ldots V_l \;$ sont des germes de
champs de vecteurs en $\;0\;$ qui engendrent Ann$(f)\;$ au voisinage de $\;0$, le lieu des z{\'e}ros communs {\`a} $\;V_1 \ldots V_l\;$ est r{\'e}duit {\`a} $\;0\;$;  r{\'e}ciproquement ceci implique la condition $\;HI$ b).

\begin{defn}
On d{\'e}finit alors $\;{\cal J}\;$ comme l'id{\'e}al {\`a} gauche de
$\;{\cal D}\;$ engendr{\'e} au voisinage de $\;0$ par
 $\;\widehat{J(f)}\;$ et Ann$(f)$.

\noindent On pose alors $\;{\cal M} = {\cal D}\big/{\cal J}$.
\end{defn}

\begin{lemma}\label{DR.1}

\noindent  On consid{\`e}re {\`a} l'origine de $\;\mathbb{C}^{n+1}\;$ des germes $\;(g_{\lambda})_{\lambda \in [1,L]} \;$ de fontions holomorphes et $\;(V_{\mu})_{\mu \in [1,M]}\;$
des germes de champs de vecteurs holomorphes. Soit $\;{\cal J}\;$ l'id{\'e}al {\`a} gauche de $\;{\cal D}_{\mathbb{C}^{n+1}}\;$ engendr{\'e}
au voisinage de $\;0\;$ par les $\;g_{\lambda}\;$ et les $\;V_{\mu}$.

\noindent Posons $\quad\widetilde{V_{\mu}} = V_{\mu} + div (V_{\mu}) \quad$ pour $\quad \mu
\in [1,M]\;$ \ o{\`u}
$$\quad  div(V):=\displaystyle  \sum^n_{i=0} \ \frac{\partial a_{i}}{\partial x_{i}}
\quad\quad{\rm si} \quad  \;V= \displaystyle  \sum^n_{i=0} a_{i} \frac{\partial}{\partial x_{i}}\ $$
 et soit $\; {\cal M} := {\cal D}_{\mathbb{C}^{n+1}}\big/{\displaystyle  {\cal J}}.$

\noindent Alors $\;DR^{n+1} ({\cal M})$, le $\;(n+1)-$i{\`e}me faisceau de cohomologie du
complexe de De Rham de $\;{\cal M}$, est isomorphe comme $\;\underline{\mathbb{C}}_X-$module au quotient 
$${\cal O}\displaystyle \Big/ {\displaystyle {\sum^{L}_{\lambda =1}}} {\cal O}. g_{\lambda} + \sum^M_{\mu=1} \widetilde{V}_{\mu} ({\cal O}).$$
 En particulier pour $\;{\cal M}\;$ holonome, le $\;\mathbb{C}-$espace vectoriel 
$${\cal O}_0 \Big/\sum^L_{\lambda=1} \; {\cal O}_0. g_{\lambda} + \sum^M_{\mu=1}
\widetilde{V}_{\mu} ({\cal O}_0)$$
sera de dimension finie, d'apr\`es [K.75].
\end{lemma}

\bigskip

\noindent \textit{\underline{Preuve}.} Comme r\'ef\'erence sur les \ $\cal{D}-$modules, le lecteur pourra consulter [Bj.93].  

\noindent Pour $\;h \in {\cal O}\;$ et $\;V\;$ un champ vecteur holomorphe  on a $\;\widetilde{V}
(h) = \widetilde{V} . h-h.V\;.$ Comme on a $\; \widetilde{V} = \displaystyle  \sum \frac{\partial}
{\partial x_{i}} \cdot a_{i}\;$ si $\; V = \displaystyle  \sum a_{i}. \frac{\partial}{\partial x_{i}},\;$ on aura, quand \ $V$ \ annule \ $f$ :
$$\widetilde{V}_{\mu} (h) \in \sum^n_{i=0} \frac{\partial}{\partial x_{i}} {\cal D}
+ {\cal J}.$$

\noindent Mais, par d{\'e}finition, 
$$DR^{n+1} ({\cal M}):= {\cal M}\Big/{\displaystyle  \sum^n_{i=0} \frac{\partial}{\partial
x_{i}} {\cal M}} \simeq {\cal D} \Big/{\displaystyle  \sum^n_{i=0} \frac{\partial}{\partial x_{i}}
{\cal D} + {\cal J}}.$$

\noindent On en d{\'e}duit que $\; \forall h \in {\cal O}\;$ et $\; \forall \mu \in [1,M],\
\widetilde{V}_{\mu}(h) \;$ est nul dans $\; DR^{n+1}( {\cal M})$.

\noindent Par ailleurs on a $\;{\cal D}\Big/{\displaystyle  \sum^{n}_{0} \frac{\partial}{\partial
x_{i}} {\cal D}} \overset{\sim}{\longrightarrow} {\cal O}\;$ et donc 
$DR^{n+1} ({\cal M}) $ \ est isomorphe \`a \ ${\cal O}\Big/{\displaystyle  {\cal O} \cap \Big({\cal J} +
\sum^{n}_{i=0} \frac{\partial}{\partial x_{i}} {\cal D}\Big)}.$

\noindent Mais on vient de voir que l'on a
$$\sum^L_{\lambda =1} {\cal O} g_{\lambda} +\sum^M_{\mu=1} \widetilde{V}_{\mu}
({\cal O})\  \subset \ {\cal O} \cap \Big({\cal J} + \sum^n_{i=0} \frac{\partial}
{\partial x_{i}} {\cal D} \Big) . $$

\noindent Il nous reste {\`a} montrer que \ $ {\cal O} \cap \Big( {\cal J} + \sum \frac{\partial}{\partial x_{i}}{\cal D} \Big)  \subset \ 
 \sum^L_{\lambda =1} {\cal O} g_{\lambda} + \sum^M_{\mu=1} \widetilde{V}_{\mu}
({\cal O}).$

\noindent Pour  \ $  \varphi \in {\cal O} \cap \Big( {\cal J} + \sum \frac{\partial}{\partial x_{i}} {\cal D} \Big) $ \ posons
\begin{equation*}
\varphi = \sum^L_{\lambda=1} P_{\alpha} g_{\lambda} + \sum^M_{\mu=1} Q_r V_{\mu}
+ \sum^n_{i=0} \frac{\partial}{\partial x_{i}} R_{i} 
\end{equation*}
o\`u \ $ P_{\lambda}, Q_{\mu} \; {\rm et} \; R_{i} \in {\cal D}.$

\noindent Ecrivons les $\; \displaystyle  \frac{\partial}{\partial x_{i}} \;$ "{\`a} gauche". On trouve
alors, grace {\`a} l'unicit{\'e} de l'{\'e}criture {\`a} gauche,
$$\varphi = \sum^L_{\lambda=1} p_{\lambda} .g_{\lambda} + \biggl[\sum^M_{\mu=1}
q_{\mu} V_{\mu}+ \sum^n_{i=0} \frac{\partial}{\partial x_{i}} r_{i}\biggr]_0$$

\noindent o{\`u} $\;p_{\lambda}, \; q_{\mu}\;$ et $\;r_{i}\;$ sont dans $\;{\cal O}\;$ et
o{\`u} le symbole $\;[\pi]_0\;$ d{\'e}signe le terme de degr{\'e} $\;0\;$ de l'op{\'e}rateur diff{\'e}rentiel
$\;\pi\;$ (avec les $\;\displaystyle  \frac{\partial}{\partial x_{i}}\;$ {\`a} gauche).

\noindent Mais on a $\;q_{\mu} V_{\mu} = \widetilde{V}_{\mu} q_{\mu} - \widetilde{V}_{\mu}
(q_{\mu})\;$ avec $\; \widetilde{V}_{\mu} \in \displaystyle  \sum^n_{i=0} \frac{\partial}
{\partial x_{i}} {\cal D}.\;$ On en conclut que 

\noindent $ \varphi = \sum^L_{\lambda=1} p_{\lambda}. g_{\lambda} - \sum^M_{\mu=1} \widetilde{V}_{\mu}
(q_{\mu}). \hfill  \blacksquare$

\begin{lemma}\label{holonome}
\noindent On se place dans la situation du lemme pr\'ec\'edent. On suppose 
\begin{enumerate}[i)]
\item que l'ensemble $\;\Big\{x \in \mathbb{C}^{n+1} \big/{\displaystyle  g_{\lambda} (x)} = 0 \quad \forall \lambda \in
      [1,L] \Big\} \;$ est un germe de courbe en $\;0.\;$ Notons le~$\;S$.
\item que pour tout point $\;y \in S\setminus \{0\}\;$ il existe $\;m_y \in [1,M]\;$ tel que
      $\;V_{m_{y}}\;$ ne s'annule pas en $\;y$.
\end{enumerate}
Alors $\;{\cal M} := {\cal D} \big/ {\cal J}\;$ est holonome.
\end{lemma}

\bigskip

\noindent \textit{\underline{Preuve}.} Il est clair que $\;$ Supp$\,{\cal M} \subset S.\;$ Au point g{\'e}n{\'e}rique $\;y\;$ de $\;S$, la vari\'et\'e caract\'eristique  ch$({\cal M})\;$ de \ $\mathcal{M}$ \ est contenue dans un fibr{\'e} en hyperplan sur $\;S\;$ (d{\'e}fini
par l'annulation du symbole principal de $\;V_{m_{y}}).\;$ On a donc $\;$ dim
ch$({\cal M}) \le n+1\;$ au dessus de \ $S \setminus\{0\}$ \  donc partout. $\hfill \blacksquare$

\subsection{Finitude de \ $Ker\,b$.}

\noindent  D{\'e}finissons maintenant les faisceaux \ $\mathcal{E}' = \Omega^{n+1}_X\big/d(Ker\,df)^n \simeq \mathcal{H}^{n+1}, \\ {\cal E}'' := \Omega^{n+1}_{X}/df_{\wedge}d \Omega^{n-1}_{X} $ \  et \ les \ espaces \ vectoriels
$$ E':= H^0_{\{0\}} (\mathcal{H}^{n+1})\quad , \quad E'' := H^0_{\{0\}} ( {\cal E}'').$$

 \noindent Soit $\;\widetilde{j} : {\cal E}'' \longrightarrow \mathcal{H}^{n+1}$ \ le
quotient {\'e}vident et notons  $\;j : E'' \longrightarrow E'\;$ l'application
qui s'en d{\'e}duit.

\begin{lemma}\label{mu}
\noindent La fl\`eche \ $ \tilde{b} : \mathcal{H}^{n+1} \rightarrow  \mathcal{E}'' $ \ donn\'ee par \ $ \tilde{b}[d\xi] = [df\wedge \xi ] $ \ est bien d\'efinie; c'est une injection  \ $\underline{\mathbb{C}}-$lin\'eaire  de conoyau \ $\Omega^{n+1}\big/df\wedge \Omega^n \simeq  \mathcal{O} \big/ J(f) \  .$
\end{lemma}

\bigskip

\noindent \textit{\underline{Preuve}.} Commen{\c c}ons par montrer que cette fl\`eche est bien d\'efinie. Consid\'erons donc un  \ $ \xi \in \Omega^{n} $ \ tel que \ $d \xi = 0$.  Alors le lemme de De Rham holomorphe permet d'\'ecrire \ $ \xi = d\eta $ \ avec \ $ \eta \in \Omega^{n-1}$. On aura donc \ $ [df\wedge \xi] = [df \wedge d\eta] = 0 $ \ dans \ ${\cal E}''.$ 

\noindent Montrons l'injectivit\'e. Si on a \ $ [df \wedge \xi] = 0 $, cela signifie qu'il existe \ $ \zeta \in \Omega^{n-1}$ \ v\'erifiant \ $ df\wedge \xi = df\wedge d\zeta \ .$ On aura donc \ $ \xi -d\zeta \in Ker\, df^n $ \ d'o\`u \ $ [d\xi] = [d(\xi -d\zeta)] = 0 $ \ dans \ $ \mathcal{H}^{n+1}$. 

\noindent L'assertion sur le conoyau de \ $\tilde{b}$ \ est \'evidente. $ \hfill \blacksquare $

\begin{lemma}\label{supp}
\noindent Sous l'hypoth{\`e}se \ $ (HI)$ \  on a au voisinage de tout point \ 
$p$ \ de \  $X \setminus \{0\}\;$:
\begin{itemize}
\item{1)} $\;d(Ker\,df)^n= \Omega^{n+1}_X $
\item{2)} $\; \Omega^n_X = d \Omega^{n-1}_X + {\rm Ker}df^n \quad $  et donc
\item{3)} $\quad df_{\wedge} d \Omega^{n-1}_X = df_{\wedge} \Omega^n_X \ .$
\end{itemize}
\end{lemma}

\noindent \textit{\underline{Preuve}.} Il suffit \'evidemment de traiter le cas o\`u \ $p\in S \setminus \{0\}.$ Soit $\;V=   \sum^n_{i=0} a_{i}. \frac{\partial}{\partial x_{i}} \;$ un champ de vecteur holomorphe au voisinage de $\;p\;$ v{\'e}rifiant $\; V.f \equiv 0\;$
et $\;V(p) \neq 0$. 

\noindent Posons $\; \alpha =   \sum^n_{i=0} a_{i} \; dx^i\;$ o{\`u} $\; dx^i = (-1)^i.d{x_{0}} {}_{\wedge} \ldots {}_{\wedge} \widehat{d{x_{i}}}
{}_{\wedge} \ldots {}_{\wedge} dx_n$. \\
 Alors $\; \alpha \in \Omega^n_X \  ,  \alpha(p) \neq
0\;$ et $\;df_{\wedge} \alpha \equiv 0$.

\noindent Pour prouver 1) on veut montrer que $\; \forall g \in {\cal O}_p\;$ il existe
$\; h \in {\cal O}_p\;$v\'erifiant \\ $\;d(h \alpha) = g \; d{x_{0}} {}_{\wedge} \ldots {}_{\wedge}
d{x_{n}} \;$ c'est-{\`a}-dire que
$$\sum^n_{j=0} \alpha_{j}. \frac{\partial h}{\partial x_{j}} + \biggl( \sum^n_{j=0}
\frac{\partial \alpha_{j}}{\partial x_{j}} \biggr).h = g.$$

\noindent Mais puisque $\;V(p) \neq 0\;$ on peut, dans un syst{\`e}me de coordonn{\'e}es
locales convenablement centr{\'e}es en $\;p,\;$ se ramener {\`a} r{\'e}soudre, $\; \forall g \in
{\cal O}_0,\;$ une {\'e}quation du type \ $\frac{\partial}{\partial x_0} h + \xi . h = g \qquad \hbox{o{\`u}} \quad h \in {\cal O}_0$ \ est inconnue et o{\`u} $\; \xi \in {\cal O}_0\;$ est donn{\'e}e. \\
Ceci est {\'e}l{\'e}mentaire.

\noindent Montrons 2). Soit $\; w \in \Omega^n_p.\;$ Alors d'apr{\`e}s 1) on peut trouver $\; h \in 
{\cal O}_0\;$ telle que $\; dw = d(h \alpha).\;$ D'apr{\`e}s De Rham on aura donc
$\; u \in \Omega^n_0 \;$ telle que $\;w-h \alpha = du.\;$ On en d{\'e}duit alors 2)
puis 3). $\hfill \blacksquare$

\bigskip

\noindent  Les faisceaux \ $\mathcal{H}^{n+1}$ \ et \ $ \tilde{b}.\mathcal{H}^{n+1} $ \ sont donc \`a support l'origine.  On pourra donc confondre le faisceau  \ $\mathcal{H}^{n+1}$ \  et l'espace vectoriel de ses sections \`a support l'origine \ $E'$. On prendra garde que le faisceau  \ ${\cal E}''$ \ n'est pas en g\'en\'eral  \`a support l'origine sous l'hypoth\`ese \ $(HI)$.

\noindent On a une suite exacte courte gr\^ace au lemme \ref{mu}
$$ 0 \rightarrow \mathcal{H}^{n+1} \overset{\tilde{b}}{\longrightarrow} {\cal E}'' \longrightarrow \Omega^{n+1}\big/ df\wedge \Omega^{n} \rightarrow 0 $$
qui donne une suite exacte courte de cohomologie \`a support l'origine, gr\^ace au lemme \ref{supp}:
$$0 \rightarrow E' \overset{\tilde{b}}{\longrightarrow}  E'' \longrightarrow H^{0}_{\lbrace 0 \rbrace}(\Omega^{n+1} \big/ df\wedge  \Omega^{n} ) \rightarrow 0. $$
Mais l'endomorphisme \ $ b : E' \rightarrow E' $ \ est la compos\'ee de \ $\tilde{b} $ \ et de \ $j$.
La premi\`ere assertion de finitude du th\'eor\`eme \ref{finitude}  est donn\'ee par la proposition suivante,
compte tenu du lemme \ref{mu} et de ce qui pr\'ec\`ede.

\begin{prop}\label{Annulation}

\noindent Sous l'hypoth\`ese \ $(HI)$ \ le noyau \ $\displaystyle H^0_{\{0\}}(X, \frac{d(Ker\,df)^n}{df\wedge d\Omega^{n-1}})$ \ de l'application  \ $ j : E'' \rightarrow E'$ \ est de dimension finie.

\noindent La dimension de \ $Ker\, j $ \ et donc \`a fortiori celle de \ $Ker\, b$ \ est major\'ee par la dimension de \ $ H^1_{\{0\}}(S, \mathcal{H}^n\big/b\mathcal{H}^n)$\footnote{Les faisceaux \ $\mathcal{H}^p$ \ ont \'et\'e introduits \`a la proposition \ref{fonction}}.
\end{prop}

\bigskip

\noindent \textit{\underline{Preuve}.} Commen{\c c}ons par montrer que l'on une suite exacte de faisceaux \`a supports dans  \ $S$.
$$0 \longrightarrow \frac{{\rm (Ker\, }df)^n \cap {\rm Ker\, }d}{{\rm (Im\, }df)^n \cap {\rm Ker\, }d}
\longrightarrow \frac{{\rm (Ker\, } df)^n}{{\rm (Im\, }df)^n} \overset{d}{\longrightarrow} 
\frac{d({\rm Ker\, }df)^n}{df_{\wedge} d\Omega^{n-1}_X} \longrightarrow 0 \ .$$
L'exactitude r{\'e}sulte simplement du fait que si
$\;w \in {\rm (Ker\, }df)^n\;$ v{\'e}rifie $\;d w = df_{\wedge}d \mu,\;$ alors $\;d(w+
df_{\wedge} \mu)=0\;$ et donc $\;[w] = [w +df_{\wedge} \mu]\;$ dans $\; \displaystyle 
\frac{{\rm (Ker\, }df)^n}{{\rm (Im\, }df)^n}\;$ est bien l'image de $\;w + df_{\wedge} \mu
\in {\rm (Ker\, }df)^n \cap {\rm Ker\, }d$. 

\smallskip

\noindent Nous aurons besoin du lemme suivant  et de son corollaire.

\begin{lemma}
\noindent  Sous l'hypoth\`ese \ $(HI)$ \ on a pour tout $\;p \in [1,n]\;$ 
$${\cal A} (p) = \left \{ 
\begin{matrix}
&{\rm (Ker\, }df)^{p-1} = {\rm (Im\, }df)^{p-1} \\
& \underline{H}^{i}_S ({\rm (Im\, }df)^{p-1})=0 \quad {\rm pour} \quad i \le n-(p-1).
\end{matrix}
\right. $$
\end{lemma}

\bigskip

\noindent \textit{\underline{Preuve}.} Pour $\;p=1,\; {\cal A} (1)\;$ est vraie car $\; {\rm (Ker\, }df)^0 = {\rm (Im\, }df)^0
=0$.

\noindent Supposons $\;{\cal A}(p)\;$ vraie avec $\;p \in [1,n-1]\;$ et montrons-la pour $\;p+1\;$ :
la suite exacte de faisceaux sur $\;X$ :

\vskip -4mm
\begin{equation*}
0 \longrightarrow {\rm (Ker\, }df)^{p-1} \longrightarrow \Omega^{p-1}_X 
\overset{\wedge df}{\longrightarrow} {\rm (Im\, }df)^p \longrightarrow 0 \tag{*}
\end{equation*}
\noindent donne pour $\;i +1 \leq n,\;$ puisque $\;\underline{H}^k_S (\Omega_X)=0\;$ pour
$\;k<n\;$ (= codim$_X S$) :
$$\underline{H}^{i}_S ({\rm (Im\, }df)^p) \subset  \underline{H}^{i+1}_S ({\rm (Im\, }df)^{p-1})$$

\noindent en utilisant $\;{\cal A}(p)\;$ qui donne $\;{\rm (Ker\, }df)^{p-1} = {\rm (Im\, }df)^{p-1}.\;$
Pour $\;i \le n-p\;$ on a $\;i+1 \le n-(p-1) \leq n$ \  et donc $\;  \underline{H}^{i+1}_S ({\rm (Im\, }df)^{p-1}) =0\;$ gr{\^a}ce a $\;{\cal A}(p)$.

\noindent Montrons maintenant que $\;{\rm (Ker\, }df)^p = {\rm (Im\, }df)^p.\;$ Comme on a d{\'e}j{\`a} obtenu que $\;\underline{H}^1_S ({\rm (Im\, }df)^p)=0,\;$ car \ $p \leq n-1$, on aura
$${\rm (Im\, }df)^p = j_*j^* \;{\rm (Im\, }df)^p$$

\noindent o{\`u} $\;j:X-S \hookrightarrow X\;$ est l'inclusion. Mais on a
 $\;j^*({\rm (Ker\, }df)^p)=j^*({\rm (Im\, }df)^p)\;$
 puisque $\;df \neq 0\;$ sur $\;X-S.\;$ De plus pour$\;n \ge2, \;S\;$ est de  codimension
$ \ge 2\;$ dans $\;X\;$ et on a $\;j_* j^* \;{\rm (Ker\, }df)^p={\rm (Ker\, }df)^p\;$ par
Hartogs\footnote{Pour $\;n=1 \  {\rm l'assertion} \,{\cal A}(1)\;$ est vraie car  vide.}.\\
 On en d{\'e}duit que $\;{\cal A} (p+1)\;$ est vraie. $\hfill \blacksquare $

\begin{cor}
\noindent Sous l'hypoth\`ese \ $(HI)$ \ on a $\; \displaystyle  H^0_{\{0\}} \big(S, \frac{{\rm Ker\, }df^n}{{\rm Im\, }df^n} \big) =0$.
\end{cor}

\bigskip

\noindent \textit{\underline{Preuve}.} Consid{\'e}rons {\`a} nouveau la suite exacte $(^*);$ comme on a 
 $$\quad H^{i}_{\{0\}}(X, \Omega^{p-1}_X)=0 \quad  \forall p \ge 1\quad  \forall i \le n$$
 on aura pour tout $\;p \le n-1\quad\;$et tout \ $i$ \ tel que $\; i+p \le n+1\;$ :
\begin{equation*}
H^{i}_{\{0\}} (X, {\rm Im\, }df^{p})=0 \ . \tag{@}
\end{equation*}
\noindent En effet c'est clair pour \ $p=0$  ; si c'est vrai pour $\;p \le n-1\ $, alors la suite
exacte (*) et $\;{\cal A} (p+1)\;$ donne l'inclusion :
$$H^{i}_{\{0\}} (X, {\rm (Im\, }df)^{p+1}) \subset  H^{i+1}_{\{0\}} (X, {\rm (Im\, }df)^{p})$$

\noindent pour $\;i+1+p \leq n+1\;$.

\noindent On en conclut, en prenant dans (@)  \ $p=n-1$ \  et \ $ i=2$ \ que, grace \`a $\;{\cal A} (n)$, $\;H^2_{\{0\}} (X, {\rm (Ker\, }df)^{n-1})=0$. La suite exacte pour $\;(^*)\;$ avec $\;p=n\;$ donne alors
$$H^1_{\{0\}} (X, {\rm (Im\, }df)^n)=0.$$

\noindent La suite exacte
$$0 \longrightarrow {\rm Im\, }df^n \longrightarrow {\rm (Ker\, }df)^n \longrightarrow
\frac{(Ker\, df)^n}{(Im\, df)^n} \longrightarrow 0$$
donne enfin, puisque $\;H^0_{\{0\}}(X, {\rm (Ker\, }df)^n)=0 \ , \ H^0_{\{0\}} \big(X, \frac{{\rm (Ker\, }df)^n}{{\rm (Im\, }df)^n}\big)=0 \ .  \hfill  \blacksquare $

\bigskip

\noindent \textit{\underline{Fin de la preuve de la proposition} \ref{Annulation}.} Remarquons maintenant que l'hypoth{\`e}se $\;(HI)\;$ implique que le faisceau \ $\mathcal{H}^n$ \ est localement constant sur \ $S^*$.  en effet pr{\`e}s de $\;y \in S^*\;$ l'existence d'un champ de vecteur holomorphe non nul en $\;y\;$ et annulant $\;f\;$ permet de choisir un syst{\`e}me de coordonn{\'e}es o{\`u} $\;f\;$ ne d{\'e}pend pas de la variable $\;x_0.\;$. Notons \ $ \pi : \mathbb{C}^{n+1} \longrightarrow \mathbb{C}^n$ \ la projection dans un tel  syst\`eme de coordonn\'ees. Comme $\;f\;$ est {\`a} singularit{\'e} isol{\'e}e dans $\; \mathbb{C}^n,\;$ on en d{\'e}duit que
$${\rm (Ker\, }df)^n \cap {\rm Ker\, }d \simeq   \pi^{-1} (\Omega^n_{\mathbb{C}^{n}}) \quad {\rm modulo} \quad df\wedge d\Omega^{n-2} \simeq d({\rm Ker\,} df)^{n-1}$$
d'o\`u notre assertion. Alors  le faisceau 
$$\displaystyle \mathcal{H}^n \big/b \mathcal{H}^n : = \frac{{\rm (Ker\, }df)^n \cap {\rm Ker\, }d}{{\rm (Im\, }df)^n \cap {\rm Ker\, }d}\;$$
 est un syst{\`e}me local sur $\;S^*$, de rang $\;\mu_{tr}$ \  pr{\`e}s de $\;y$, o\`u \ $\mu_{tr}={\rm dim}_{\mathbb{C}} \big(\Omega^n_{\mathbb{C}^{n}}\big/ df_{\wedge} \Omega^{n-1}_{\mathbb{C}^{n}})$  \ est le rang sur \ $\mathbb{C}[[b]]$ \ du faisceau localement constant sur \ $S^*$ \ de (a,b)-modules associ\'e au faisceau localement constant de pr\'e-(a,b)-modules (sans torsion) \ $\mathcal{H}^n $.

\noindent Comme $\;S\;$ est une r{\'e}union finie de disques topologiques ayant m{\^e}me centre
(et sinon disjoints) on en d{\'e}duit que $\; H^1_{\lbrace 0 \rbrace}(S, \mathcal{H}^n \big/b \mathcal{H}^n)$ \ est de dimension finie sur \ $ \mathbb{C}  .$

\noindent La suite exacte de cohomologie \`a support l'origine de la suite exacte 
$$0 \longrightarrow \frac{{\rm (Ker\, }df)^n \cap {\rm Ker\, }d}{{\rm (Im\, }df)^n \cap {\rm Ker\, }d}
\longrightarrow \frac{{\rm (Ker\, } df)^n}{{\rm(Im\, }df)^n} \overset{d}{\longrightarrow} 
\frac{d{\rm (Ker\, }df)^n}{df_{\wedge} d\Omega^{n-1}_X} \longrightarrow 0 \ .$$
donne alors le "tron{\c c}on" exact :
$$ 0 \simeq H^0_{\lbrace 0 \rbrace}\big(S,\frac{{\rm (Ker\, } df)^n}{{\rm (Im\, }df)^n}\big) \rightarrow H^0_{\lbrace 0 \rbrace}\big(S, \frac{d{\rm (Ker\, }df)^n}{df_{\wedge} d\Omega^{n-1}_X}\big) \rightarrow  H^1_{\lbrace 0 \rbrace}\big(S, \mathcal{H}^n \big/b \mathcal{H}^n \big) $$
qui permet de conclure \`a la finitude de l'espace vectoriel
$\displaystyle H^0_{\lbrace 0 \rbrace}\big(S, \frac{d{\rm (Ker\, }df)^n}{df_{\wedge} d\Omega^{n-1}_X}\big)$ \  et donc \`a la finitude de \ $ Ker\, j . \hfill \blacksquare $

\bigskip

\noindent  \textit{\underline{Fin de la d\'emonstration du th\'eor\`eme \ref{finitude}}.} Il nous reste \`a montrer que le conoyau de \ $ b : E' \rightarrow E' $ \ est de dimension finie. Ceci va r\'esulter de la finitude du quotient \ $ E'' \big/ \tilde{b}E' $ \ qui a d\'eja \'et\'e obtenue, et  de la finitude donn\'ee par les lemmes \ref{DR.1}, \ref{holonome} et le lemme \ref{DR.2} ci-dessous. La suite exacte suivante  permet alors de conclure
$$ E''\big/\tilde{b}E' \overset{j}{\longrightarrow} E'\big/bE' \longrightarrow E'\big/ jE'' + bE'  \rightarrow 0  \ . $$
puisque que l'on a d\'eja obtenu la finitude de \ $ Ker\, j $.

\begin{lemma}\label{DR.2}
\noindent  Sous l'hypoth\`ese \ $(HI)$ \ on a un isomorphisme d'espace vectoriel
$$ DR^{n+1}(\mathcal{M}) \tilde{\longrightarrow} E'\big/ jE'' + bE' \simeq E'\big/ jE'' . $$
\end{lemma}

\bigskip

\noindent \textit{\underline{Preuve}.} Rappelons que  le $\;{\cal D}-$module $\;{\cal M}\;$ est d{\'e}fini comme le quotient 
$\;{\cal D} \big/ {\cal J}\;$ o{\`u} $\;{\cal J}\;$ est l'id{\'e}al {\`a} gauche de $\;{\cal D}\;$
engendr{\'e} par $\;\widehat{J(f)}\;$ et Ann$f \subset T_X\; .$ D'apr\`es les lemmes  \ref{DR.1} et \ref{holonome}, il  est holonome  et  on a  l'isomorphisme
$$DR^{n+1} ({\cal M}) \overset{\sim}{\longrightarrow} {\cal O}\Big/{\displaystyle 
\sum^M_{j=1} {\cal O} g_{j}} + \sum^l_{j=1} \widetilde{V}_{j} ({\cal O})$$

\noindent o{\`u} $\;g_1 \ldots g_M\;$ engendrent $\;\widehat{J(f)} \;$ sur $\;{\cal O},\;$o\`u 
$\;V_1 \ldots V_l\;$ engendrent Ann$(f)\;$ sur $\;{\cal O}_X\;$ et o\`u on associe au champ de vecteur \ $V$ \ l'op\'erateur diff\'erentiel \ $\tilde{V}$ \ d'ordre 1 d\'efini par 
$$\widetilde{V}(h) = V . h +  div (V).h$$

\noindent Montrons que ceci est bien lin{\'e}airement
isomorphe {\`a} $\;E'\big/\displaystyle  j(E'')+bE' \simeq E'\big/ jE''$  \ comme $\;\mathbb{C}-$espace vectoriel.

\smallskip

\noindent D'abord $\;w \in \Omega^{n+1}_{X,0}\;$ donne une classe dans $\;E'' =H^0_{\{0\}}
( \Omega^{n+1} \big/ {\displaystyle  df_{\wedge}d \Omega^{n-1}_{X}})\;$ si et 
seulement si localement dans  \  $  X- \lbrace 0 \rbrace, \ w$ \  est dans $\;df_{\wedge} d\Omega^n_X\;$.
D'apr{\`e}s l'{\'e}galit{\'e} $\;df_{\wedge} d \Omega^{n-1}_X = df_{\wedge} \Omega^n_X\;$ sur \ 
$X^*$ \  prouv{\'e}e au lemme \ref{supp}, \ $w$ \ induira une classe de \ $E''$ \  si et seulement si $\;w \in \widehat{J(f)}.\Omega^{n+1}_X$. 

\noindent Si $\;v_1 \ldots v_l \in {\rm (Ker}\, df)^n\;$ correspondent {\`a} $\;V_1 \ldots V_l,\;$ on aura \\
$\sum^l_{j=1} d({\cal O} v_{j}) = \biggl(\sum^l_{j=1} \widetilde{V}_{j} ({\cal O})
\biggr).dx.$\\
 On en d{\'e}duit, puisque la relation \ $j_{\circ}\tilde{b} = b$ \ montre que \ $ bE' \subset jE''$, les isomorphismes 
\begin{equation*}
E'\Big/ {  j E'' + bE'}  \simeq \Omega^{n+1}_X \Big/{ \widehat{J(f)}.
\Omega^{n+1}_X+ \biggl(\sum^l_{j=1} \widetilde{V}_{j} ({\cal O})\biggr).dx}  \simeq {\cal O} \Big/ { \widehat{J(f)} + \sum^l_{j=1} \widetilde{V}_{j}
({\cal O}) }. 
\end{equation*} 
 Ceci prouve notre assertion. \ $\hfill \blacksquare $

\subsection{G\'en\'eralisation de la formule de J. Milnor.}

Une cons\'equence simple de ce qui pr\'ec\`ede est la g\'en\'eralisation suivante de la formule de J. Milnor, donnant, sous nos hypot\`eses, la dimension du $n-$i\`eme groupe de cohomologie de la fibre de Milnor de \ $f$ \ \`a l'origine.

\begin{cor}\label{Miln.gener.}

Sous l'hypoth\`ese (HI) l'espace vectoriel \ $ E' : = \mathcal{H}^{n+1}_0 \simeq H^0_{\lbrace 0 \rbrace}(Y, \mathcal{H}^{n+1})$ \ est un pr\'e-(a,b)-module de rang \ $r$ \ v\'erifiant
$$  \dim E'/b.E' =  \mu(f) + \nu(f) - \gamma + \delta \quad {\rm et} \quad \dim H^n(F_0, \mathbb{C}) = r =  \mu(f) + \nu(f) - \gamma.$$
o\`u \ $F_0$ \ d\'esigne la fibre de Milnor de \ $f$ \ \`a l'origine et o\`u \ $\gamma$ \ et \ $\delta$ \ sont les dimensions respectives de \ $Ker\, j$ \ et de \ $Ker\, b$.
%En particulier, sous les hypoth\`eses \ $(HI)$ \ {\bf et}  la condition \ $(P)$ \   on aura
%$$\dim H^n(F_0, \mathbb{C}) = r =  \mu(f) + \nu(f). $$
\end{cor}

\bigskip

\noindent \textit{D\'emonstration.} Rappelons que le rang d'un pr\'e-(a,b)-module est, par d\'efinition, le rang du (a,b)-module associ\'e. C'est donc le rang du fibr\'e de Gauss-Manin de \ $f$ \ \`a l'origine et  ce rang est \'egal \`a la dimension de \ $H^n(F_0, \mathbb{C})$.
 La suite exacte
$$ 0 \to Ker\, j \to E'' \overset{j}{\to} E' \to E'/jE'' \to 0 $$
donne la suite exacte d'espaces vectoriels de dimensions finies, gr\^ace au th\'eor\`eme de finitude \ref{finitude}
$$ 0 \to Ker\, j/(Ker\, j \cap \tilde{b}.E') \to E''/ \tilde{b}E'  \to E'/b.E' \to E'/jE'' \to 0 .$$
On a \ $ Ker\, b \simeq Ker\, j \cap \tilde{b}.E' $ \ via \ $\tilde{b}$ \ puisque \ $\tilde{b}$ \ est injective et que  \ $b = j_{\circ} \tilde{b}$ \ dans \ $E'$. On obtient alors, puisque \ $ \dim E''/\tilde{b}.E' = \mu(f), \dim E'/jE'' = \nu(f) $, la relation : 
$$ \dim (E'/b.E') + \gamma - \delta = \mu(f) + \nu(f).$$
La suite exacte \ $ 0 \to B(E') \to E' \to E'/B(E') \to 0 $ \ donne alors la relation
$$ \dim E'/b.E' = \dim \big(E'/B(E') + b.E' \big) + \dim \big(B(E')/b.B(E')\big).$$
Comme on a \ $ r = \dim \big(E'/B(E') + b.E'\big) $ \ et \ $ \dim B(E')/b.B(E') = \dim Ker\, b $, on conclut facilement. $\hfill \blacksquare$

\bigskip

\noindent {\bf Remarque importante.}

\bigskip

\noindent Les entiers \ $\gamma$ \ et \ $\delta$ \ sont bien g\'enants, aussi est-il interessant de disposer de situations o\`u l'on sait qu'ils sont nuls.

\noindent Pour avoir l'injectivit\'e de l'application \ $E''\overset{j}{\longrightarrow} E' $, c'est \`a dire la nullit\'e de \ $\gamma$, une condition n\'ecessaire et suffisante est donn\'ee par l'inclusion
\begin{equation}
 d(Ker\, df^n) \cap \widehat{J}(f).\Omega^{n+1} \subset  df\wedge d\Omega^{n-1}_X ,  \tag{P}
\end{equation}
alors qu'une condition n\'ecessaire et suffisante pour l'injectivit\'e de \ $ E' \overset{b}{\longrightarrow}E' $, c'est \`a dire pour la nullit\'e de \ $\delta$, est donn\'ee par l'inclusion 
\begin{equation*}
 d(Ker\, df^n) \cap (df\wedge \Omega^n_X) \subset df\wedge d\Omega^{n-1}_X .  \tag{P'}
\end{equation*}
Bien sur on a toujours \ $\delta \leq \gamma$ \ et la  condition \ $(P)$ \ implique la condition  \ $(P')$.

\begin{cor} Sous les hypoth\`eses \ $(HI)$ \ et  la condition \ $(P)$ \   on aura
$$ rg(E') = {\rm dim\, } E'\big/ bE'  = \mu(f) + \nu(f), $$
o\`u \ $ \mu(f) = {\rm dim\,} \widehat{J(f)}\big/ J(f) $ \ et o\`u \ $ \nu(f) = {\rm dim\, } DR^{n+1}(\mathcal{M}). $\\
Dans ces conditions la cohomologie de degr\'e \ $n$ \ de la fibre de Milnor de \ $f$ \ \`a l'origine est de dimension \ $ \mu(f) + \nu(f) $.
\end{cor}

\noindent Nous montrerons au paragraphe  4   que cette condition est souvent v\'erifi\'ee. 

\noindent En fait je ne connais pas d'exemple de cas o\`u \ $f$ \ v\'erifie l'hypoth\`ese \ $(HI)$ \ et o\`u la  condition \ $(P)$ \ n'est pas r\'ealis\'ee.

% Remarquons enfin qu'il r\'esulte facilement de l'hypoth\`ese \ $(HI)$ \footnote{Sans m\^eme supposer que la condition \ $(P) $ \ est v\'erifi\'ee.} \ que l'on a aussi, pour chaque \ $ x \in S^*$ 
%$$ dim_{\mathbb{C}} \mathcal{H}^n\big/b\mathcal{H}^n  = dim_{\mathbb{C}} H^{n-1}(F_x) $$
% o\`u \ $F_x$ \ d\'esigne la fibre de Milnor de \ $f$ \ en  \ $x$.

\subsection{ Exemple.} Voici un  exemple simple (avec \ $ n=1$)\footnote{ Notons que $\;(HI)\;$ est toujours v{\'e}rifi{\'e}e pour $\;n=1$. C'est \'egalement le cas pour la condition \ $(P)$ \ comme on le montrera au paragraphe 4. } : 
$$\;f(X,Y) =X^3(X^3+Y^3) .$$
Il est facile de voir que $\;{\rm Ker\, }df^1\big/{\displaystyle {\cal O}.df}\;$ est engendr{\'e} 
par la $\;1-$forme associ{\'e} au champ de vecteur
$$V= XY^2 \frac{\partial}{\partial X} - (2X^3+Y^3) \frac{\partial}{\partial Y}$$

\noindent annulant $\;f\;$ et de divergence $\;{\rm div} \; V =-2Y^2$.

\noindent On a \  $J(f) = X^2.(2X^3 + Y^3 , XY^2) , \ \widehat{J(f)} = (X^2)\ $ et $\ {\rm dim}_{\mathbb{C}} \; \widehat{J(f)} \big/{\displaystyle  J(f)}=\ \mu(f) = 9 \ $. Pour calculer  ${\rm dim}_{\mathbb{C}}\;\Big(\mathbb{C} \{X,Y\}\big/{\displaystyle  (X^2)+\widetilde{V}\big(\mathbb{C} \{X,Y\}\big)}\Big)$ on  regarde l'identit{\'e}
$$\widetilde{V}(X^p Y^q)=(p-q-2) X^p Y^{q+2}-2q X^{p+3} Y^{q-1}\;;$$

\noindent on peut donc r{\'e}duire $\;X^{a} Y^b\;$ {\`a} $\;X^{a+3} Y^{b-3} \;$ pourvu que
$\; a\neq b\;$ et $\;b \ge 2$.

\noindent Si $\;a=b \ge 2 \;$ on est dans $\;(X^2)$. Il reste donc $\;1,X,Y\;$ et
$\;X^pY\;$ pour $\; p <2$.

\noindent Donc $\;1,X,Y,XY\;$ donne une base et $\;\nu(f) =4$.

\noindent Alors  $\; \mathcal{L}(E')$, le compl\'et\'e \ $b-$adique de \ $E'$ (voir le paragraphe 2),  est un $\;(a,b)-$module de rang 13. Une base de
$\;E'\big/_{\displaystyle  bE'} \;$ est donn{\'e}e par
$$1,X,Y,XY,X^2,X^3,X^4,X^5,X^2Y,X^3Y,X^4Y,X^5Y,X^2Y^2$$

\noindent et on a $\;a[1]=\displaystyle  \frac{1}{3} b[1]\;\ldots\;a[X^pY^q]= \frac{p+q+2}{6} b[X^pY^q]$.

\section{ Quelques cas sans torsion.}

\subsection{}

  Le calcul de \ $\mu(f) $ \ est en g\'en\'eral assez simple (comparable au calcul du nombre de Milnor dans le cas d'une singularit\'e isol\'ee). Le calcul de \ $ \nu(f) $ \ bien que l'on ait donn\'e un isomorphisme  de
\ $DR^{n+1}(\mathcal{M}) $ \ sur un espace vectoriel plus "concret" au lemme \ref{DR.1}, est en pratique beaucoup plus d\'elicat car l'espace vectoriel consid\'er\'e n'a pas une structure de module sur \ $ \mathcal{O}_0 $.

\bigskip

\noindent  Sans la condition \ $(P)$ \ pour \ $f$, on obtient seulement une majoration de la dimension du quotient \ $E'/b.E'$, ce qui rend  les calculs encore bien plus p\'enibles. D'o\`u l'inter\^et de savoir, pour une sous-famille aussi large que possible de la famille des fonctions v\'erifiant l'hypoth\`ese \ $(HI)$,  que la condition \ $(P)$ \ est satisfaite. Ceci est d'autant plus utile que la v\'erification de cette condition  est  d\'elicate, comme on peut s'en convaincre sur les exemples. Nous proposons donc de donner dans ce paragraphe deux r\'esultats qui permettent de voir que la condition \ $(P)$ \ est tr\`es souvent v\'erifi\'ee sous l'hypoth\`ese \ $(HI)$. 

\noindent Le premier r\'esultat (th\'eor\`eme \ref{P}) montre que la propri\'et\'e \ $(P)$ \ est toujours v\'erifi\'ee pour \ $ n=1 $, c'est \`a dire pour les  germes de fonctions holomorphes  (r\'eduites ou non) \`a l'origine  de \ $\mathbb{C}^2$. 

\noindent Le second (la proposition \ref{T.S.1}) assure que la suspension d'une fonction \`a singularit\'e isol\'ee avec une fonction qui v\'erifie \ $ (HI) \ {\rm et} \ (P) $ \ v\'erifie \'egalement \ $ (HI) \ {\rm et} \ (P).$

%%%%%%%%%%%%%%%%%%%%%%%%%%%%%%%%%%%%%%%%%%%%%%%%%%%%%%%%%%%%%%%%%%%%%%%%%%%%%%%%%%%%%%%%%%%%%%%%%%%%%%%%

\subsection{Courbes planes.}

\begin{lemma}[R\'ecurrence tordue.]\label{rec.tord.}
\noindent Soient \ $p_1, \cdots, p_k$ \ des entiers \ $\geq 2$ \ et soient  \ $\phi_1 \cdots \phi_k $ \ des fonctions strictement croissantes 
\begin{equation*}
        \qquad    \phi_j  : [0,p_{j}-1]\cap \mathbb{N}  \longrightarrow [0,1] 
\end{equation*}
v\'erifiant 
$$  \phi_{j}(p_{j}-1) = 1 ,\ \forall j \in [1,k]  . $$
Consid\'erons des propositions \ $ A(\sigma_1 ,\cdots ,\sigma_k ) $ \ ind\'ex\'ees par les entiers
 $$ \sigma_j \in [0,p_{j}-1] \cap \mathbb{N}  . $$
\noindent On suppose 
\begin{itemize}
\item{1)}  $ A(0,\cdots , 0) $ \ est vraie .
\item{2)}  l'implication \ $ A(\sigma_1 ,\cdots ,\sigma_k ) \Rightarrow  A(\sigma_1 ,\cdots ,\sigma_{j}+1, \cdots ,
\sigma_k ) $ \ est vraie si les deux conditions suivantes sont satisfaites 
\end{itemize}
\begin{align*}
\qquad  &  a)  \ \sigma_j  \leq p_{j}-2  \\
\qquad  &  b)  \ \phi_j (\sigma_j ) =\min_{l\in [1,k]} \lbrace \phi_l (\sigma_l ) \rbrace 
\end{align*}
Alors la proposition \ $ A(p_{1}-1, \cdots , p_{k}-1 ) $ \ est vraie.
\end{lemma}

\smallskip

\noindent {\bf On prendra garde que nous n'affirmons pas ici que \ $A(\sigma_1,\cdots ,\sigma_k) $ \ est vraie pour toutes les valeurs de \ $(\sigma_1,\cdots ,\sigma_k) .$}

\noindent La preuve (\'el\'ementaire) est laiss\'ee en exercice au lecteur.

\begin{thm}\label{P}
Soit  \ $f$ \ un germe non nul de fonction holomorphe \`a l'origine de \ $\mathbb{C}^{2} $, v\'erifiant \ $ f(0) = 0 $. \\
Alors \ $f$ \ v\'erifie \ $(HI)$ \ et \ $(P)$.
\end{thm}

\bigskip

\noindent Ce th\'eor\`eme est une cons\'equence imm\'ediate des deux propositions suivantes.

\begin{prop}\label{a.b.s.}
 Soit  \ $f$ \ un germe de fonction holomorphe \`a l'origine de \ $\mathbb{C}^{2} $, v\'erifiant \ $ f(0) = 0 $ \ et que l'on supposera non identiquement nul. Ecrivons 
$$ f = u_{1}^{p_1}\cdots u_{k}^{p_k}.\psi  $$
o\`u \ $ u_1 \cdots u_k $ \ sont des germes irr\'eductibles en \ $0$ \ deux \`a deux distincts, o\`u \ $\psi$ est un germe nul en \ $0$ \ suppos\'e r\'eduit (donc \`a singularit\'e isol\'ee en \ $0$ ) et o\`u les entiers \ $ p_1, \cdots , p_k $ \ sont des entiers au moins \'egaux \`a deux. On suppose de plus qu'aucun des \ $u_j $ \ ne divise \ $ \psi $.

\noindent Alors on a :
\begin{itemize}
\item{1)} \quad  $ Ker df^1  = \mathcal{O}.\alpha $ \ o\`u 
$$  \alpha : = \sum_{l=1}^{k}  p_l .u_1 \cdots \hat{u}_l \cdots u_k .\psi .du_l \  + \  u_1 \cdots u_k .d\psi  $$
\item{2)} Pour tout \ $ h \in \mathcal{O}_{\mathbb{C}^{2},0} $ \ v\'erifiant \ $  d(h.\alpha) \in (u_{1}^{p_{1}-1} \cdots u_{k}^{p_{k}-1} ).\Omega^{2}_{\mathbb{C}^{2},0} $ \  on a  $ h \in  (u_{1}^{p_{1}-1} \cdots u_{k}^{p_{k}-1} ). \mathcal{O}_{\mathbb{C}^{2},0} . $
\end{itemize}
\end{prop}

\noindent Nous allons maintenant traiter le cas o\`u \ $\psi(0) \not= 0 $ \ dans la proposition pr\'ec\'edente. On peut alors supposer que \ $ \psi \equiv 1 $.

\begin{prop}\label{s.b.s.}
Dans la m\^eme situation que la proposition pr\'ec\'edente, supposons maintenant que 
 \ $ \psi \equiv 1 $. Alors on a :
\begin{itemize}
\item{1)} \  $ Ker df^1  = \mathcal{O}.\alpha $ \ o\`u  \ $  \alpha : = \sum_{l=1}^{k}  p_l .u_1 \cdots \hat{u}_l \cdots u_k  .du_l  .$
\item{2)} \  Pour tout \ $ h \in \mathcal{O}_{\mathbb{C}^{2},0} $ \ v\'erifiant  \ $  d(h.\alpha) \in (u_{1}^{p_{1}-1} \cdots u_{k}^{p_{k}-1} ).\Omega^{2}_{\mathbb{C}^{2},0} $ \ 
on peut \'ecrire \ $ h = h_0 + h_1.u_1 ^{p_1 -1} \cdots u_k ^{p_k -1} $ \ avec  \ $ h_1 \in \mathcal{O}_{\mathbb{C}^2 ,0 } $ \ et \ $  h_0 \in \mathcal{O}_{\mathbb{C}^2 ,0 } $ \ v\'erifiant \ $ d(h_0 .\alpha ) \equiv 0 . $
\end{itemize}
\end{prop}

%%%%%%%%%%%%%%
%%%%%%%%%%%%%%

\noindent \textit{\underline{Preuve} de \ref{a.b.s.}.}  Prouvons d\'ej\`a l'assertion  1). La suite exacte de faisceaux coh\'erents 
\begin{equation*}
0 \rightarrow (Ker\,df)^1 \rightarrow \Omega^{1}_{\mathbb{C}^{2} } \rightarrow\Omega^{2}_{\mathbb{C}^{2}}\rightarrow \Omega^{2}_{\mathbb{C}^{2}}/df\wedge \Omega^{1}_{\mathbb{C}^{2}} \rightarrow 0 
\end{equation*}
montre que le faisceau \ $(Ker\, df)^1 $ \ est localement libre de rang 1  d'apr\`es le th\'eor\`eme de
Hilbert pour \ $\mathbb{C}^2  .$
\noindent Comme on a \ $  df = u_{1}^{p_{1}-1} \cdots u_{k}^{p_{k}-1}.\alpha $ on a \ $ \mathcal{O}.\alpha \hookrightarrow (Ker\, df)^1 . $ \ Il s'agit de voir que cette inclusion est une \'egalit\'e. Mais comme on sait que
\ $(Ker\, df)^1 $ \ est libre de rang 1 pr\`es de l'origine, il suffit de prouver cette \'egalit\'e en dehors de l'origine d'apr\`es Hartogs. Quand \ $df$ \ ne s'annule pas c'est clair car \ $df$ \ et \ $\alpha$ \ diff\`erent d'un facteur
inversible et on a \ $ (Ker\, df)^1 = \mathcal{O}.df $ \ pr\`es d'un tel point.

\noindent Il nous suffit donc de montrer cette \'egalit\'e pr\`es du point g\'en\'erique de \ $u_j = 0 $ \ pour chaque \ $ j \in [1,k] . $ \ Au voisinage d'un tel point \ $ df $ \ ne diff\`ere de \ $ u_{j}^{p_{j}-1}.\alpha $ \ que par un facteur inversible ; on a donc  \ $ (Ker\, df)^1 = Ker (\wedge \alpha) .$ \ Mais pr\`es d'un tel point \ $\alpha $ \ ne s'annule pas car \ $ p_j u_1 \cdots \hat{u}_j \cdots u_k .\psi.du_j $ \ est non nulle alors que la somme 
$$  \sum_{i \not= j } p_l .u_1 \cdots \hat{u}_l \cdots u_k .\psi.du_l  $$
est nulle en ce point  (puisqu'il y a \ $u_j $ \ en facteur ).

\noindent On aura donc \ $ Ker( \wedge \alpha) = \mathcal{O}.\alpha $ \ pr\`es d'un tel point, ce qui ach\`eve la preuve de l'assertion  1).

\noindent Pour prouver l'assertion  2)  nous allons utiliser la "r\'ecurrence tordue" donn\'ee au lemme pr\'ec\'edent  avec l'assertion suivante  sur \ $ h \in \mathcal{O}_{\mathbb{C}^{2},0} $ 
\begin{equation*}
A(\sigma_1 ,\cdots ,\sigma_k ) : = \Big(\lbrace d(h.\alpha) \in u_{1}^{\sigma_{1}}\cdots u_{k}^{\sigma_{k}}.\Omega^{2}_{\mathbb{C}^{2},0} \rbrace \Longrightarrow \lbrace h \in u_{1}^{\sigma_{1}}\cdots u_{k}^{\sigma_{k}}.\mathcal{O}_{\mathbb{C}^{2},0}\rbrace \Big)
\end{equation*}
Bien sur,\ $ \sigma_j \in [0,p_j -1] , \quad \forall j \in [1,k] $ \ et  \ $ A(0,\cdots ,0) $ \ est vraie. Nous utiliserons les fonctions \ $ \phi_j (x) = \frac{x+1}{p_j } \quad \rm{pour} \quad j \in [1,k] . $ 

\noindent Supposons donc que \ $ A(\sigma_1,\cdots ,\sigma_k) $ \ soit vraie et que les conditions a) et  b)
de la "r\'ecurrence tordue " soient satisfaites. Consid\'erons alors \ $ h \in \mathcal{O}_{\mathbb{C}^{2},0}$ \ tel que 
 $$  d(h.\alpha) \in (u_{1}^{\sigma_{1}} \cdots u_j ^{\sigma_j +1} \cdots  u_{k}^{\sigma_{k}} ).\Omega^{2}_{\mathbb{C}^{2},0} .  $$
Alors gr\^ace \`a \ $ A(\sigma_1,\cdots ,\sigma_k) $ \ on peut \'ecrire 
$$  h = u_{1}^{\sigma_{1}}\cdots u_{k}^{\sigma_{k}}.h'  $$
et il s'agit essentiellement de montrer que l'hypoth\`ese 
\begin{equation*}
d(u_{1}^{\sigma_{1}}\cdots u_{k}^{\sigma_{k}}.h' .\alpha)  \in (u_{1}^{\sigma_{1}} \cdots u_j ^{\sigma_j +1} \cdots  u_{k}^{\sigma_{k}} ).\Omega^{2}_{\mathbb{C}^{2},0} \qquad    \tag{*}
\end{equation*}
et les conditions  a)  et  b)  permettent de conclure que l'on a \ $ h' \in (u_j) . $ 

\noindent Quelques calculs p\'enibles montrent que l'hypoth\`ese (*) donne que la forme

\begin{align*}
\omega \  =  & \  h'. \sum_{i\not=j , i=1}^{k} (p_j .\sigma_i - p_i .\sigma_j ). u_1 \cdots \hat{u}_i \cdots \hat{u}_j \cdots u_k .\psi .du_i \\
\ & - h' .\sigma_j u_1 \cdots \hat{u}_j \cdots u_k .d\psi 
\  + p_j .u_1 \cdots \hat{u}_j \cdots u_k .\psi .dh'  \\
\  & + h' .\sum_{i\not=j , i=1}^{k} (p_j - p_i ).u_1 \cdots \hat{u}_j \cdots \hat{u}_i \cdots u_k .\psi .du_i 
\    + (p_j -1).h'.u_1 \cdots \hat{u}_j \cdots u_k .d\psi 
\end{align*}
v\'erifie \ $ \omega \wedge du_j  \in (u_j ) .\Omega^{2}_{\mathbb{C}^{2},0}  .$ 

\noindent Mais en regroupant les termes on obtient 
\begin{align*} 
\omega \ = & \ h' \sum_{i\not=j} (p_j (\sigma_i +1) - p_i (\sigma_j +1)).u_1 \cdots \hat{u}_i \cdots \hat{u}_j \cdots u_k .\psi .du_i \\
\qquad  &  + (p_j - \sigma_j -1).h'.u_1 \cdots \hat{u}_j \cdots u_k .d\psi 
\  + p_j .u_1 \cdots \hat{u}_j \cdots u_k .\psi .dh' 
\end{align*}
Comme on a suppos\'e \ $ p_j (\sigma_i +1) - p_i (\sigma_j +1) \geq 0 \quad \forall i \in [1,k] $ \ , c'est \`a dire que \ $ \phi_j (\sigma_j ) = \frac{\sigma_j +1}{p_j } $ \ est l'infimum des \ $ \frac{\sigma_i +1}{p_i} $ \ pour \ $i \in [1,k] $, la fonction
\begin{equation*}
\lambda = (h')^{p_j}.\psi ^{p_j -\sigma_j -1 }. \prod_{i\not=j} u_i ^{ p_j (\sigma_i +1) - p_i (\sigma_j +1)}
\end{equation*}
est holomorphe au voisinage de l'origine dans \ $\mathbb{C}^2 $ \ et la condition
$$\omega \wedge du_j  \in (u_j ) .\Omega^{2}_{\mathbb{C}^{2},0}  $$
donne \ $ d\lambda \wedge du_j \in (u_j ) .\Omega^{2}_{\mathbb{C}^{2},0} \ . $ 

\noindent Donc la restriction de \ $\lambda$ \  \`a la courbe lisse et connexe \ $ \lbrace u_{j} = 0 \rbrace \setminus \lbrace 0 \rbrace $ \ est constante et on a 
$$   \lambda  =  \lambda_0 + u_j .L  $$
avec \ $ \lambda_0 \in \mathbb{C}$ \   et \ $ L \in \mathcal{O}_{\mathbb{C}^2 ,0} \ . $

\noindent Comme on a \ $ \sigma_j +1 < p_j $ \ et \ $ \psi(0) = 0 $, la fonction \ $\lambda $ \ est nulle
en  0  et \ $\lambda_0 = 0  . $ \ Alors \ $ \lambda \in (u_j ) $ \ et, puisque cet id\'eal est premier et ne contient
ni \ $\psi$ \ ni \ $ u_i \  , \forall i\not=j \  , i  \in [1,k] $, on aura \ $ h' \in (u_j ) $ \ ce qui prouve  \ $A(\sigma_1 \cdots \sigma_j +1 \cdots \sigma_k ). \hfill \blacksquare $

\bigskip

\noindent On remarquera que la preuve de \ref{a.b.s.}  utilise de fa{\c c}on essentielle le fait que \ $ \psi(0) = 0 $. Le cas \ $\psi \equiv 1 $ \ est donc plus d\'elicat.

\bigskip

\noindent \textit{\underline{Preuve} de \ref{s.b.s.}.} Le point  1)  se prouve de fa\c con analogue au  1) de la proposition pr\'ec\'edente.

\noindent Pour le point  2) nous allons utiliser la m\^eme strat\'egie que pr\'ec\'edement mais avec 
\begin{align*}
A(\sigma_1 ,\cdots ,\sigma_k ) : = & \Big( \lbrace d(h.\alpha) \in u_{1}^{\sigma_{1}}\cdots u_{k}^{\sigma_{k}}.\Omega^{2}_{\mathbb{C}^{2},0} \rbrace \Longrightarrow  \\
\lbrace h = h_0 + h_1.u_1 ^{p_1 -1} \cdots u_k ^{p_k -1}\quad \rm{avec} &  \quad\    h_0 \in \mathcal{O}_{\mathbb{C}^2 ,0 }  \  \rm{v\acute{e}rifiant }\ d(h_0 .\alpha ) \equiv 0 \ \rbrace \Big) . 
\end{align*}
Comme \ $ A(0,\cdots ,0) $ \ est trivialement v\'erifi\'ee, montrons que si \ $ j\in [1,k] $ \ v\'erifie 
\begin{equation*}
  \frac {\sigma_j +1}{p_j}  = \min_{l\in [1,k]} \lbrace \frac{{\sigma_l +1}}{p_l} \rbrace < 1 
\end{equation*}
l'implication
$$  A(\sigma_1 ,\cdots , \sigma_k ) \Longrightarrow  A(\sigma_1 ,\cdots ,\sigma_j +1 ,\cdots ,\sigma_k ) $$
est vraie. 

\noindent Soit donc \ $ h \in \mathcal{O}_{\mathbb{C}^2 ,0 } $ \ tel que 
$$  d(h.\alpha) \in (u_{1}^{\sigma_{1}} \cdots u_j ^{\sigma_j +1} \cdots  u_{k}^{\sigma_{k}}) .\Omega^{2}_{\mathbb{C}^{2},0} $$
et supposons que \ $ A(\sigma_1 ,\cdots ,\sigma_k ) $ \ soit vraie. Alors on peut \'ecrire 
$$ h = h_0 + h_1.u_1 ^{\sigma_1 } \cdots u_k ^{\sigma_k } $$
avec  \ $ h_1 \in \mathcal{O}_{\mathbb{C}^2 ,0 } $ \ et \ $  h_0 \in \mathcal{O}_{\mathbb{C}^2 ,0 } $ \ v\'erifiant \ $ d(h_0 .\alpha ) \equiv 0 \  . $\ L'hypoth\`ese est alors que 
$$  d(h_1.\alpha) \in (u_{1}^{\sigma_{1}} \cdots u_j ^{\sigma_j +1}\cdots  u_{k}^{\sigma_{k}}).\Omega^{2}_{\mathbb{C}^{2},0} \  . $$
Les calculs de la proposition pr\'ec\'edente restent valables ( avec \ $ \psi \equiv 1 $) \ et conduisent \`a
$$ \omega \wedge du_j  \in  u_j .\Omega^{2}_{\mathbb{C}^{2},0} \  . $$
avec 
\begin{equation*}
\omega =   h_1 . \sum_{i\not=j} (p_j (\sigma_i +1) - p_i (\sigma_j +1)).u_1 \cdots \hat{u}_i \cdots \hat{u}_j \cdots u_k .du_i  + p_j .u_1 \cdots \hat{u}_j \cdots u_k .dh_1
\end{equation*}
Si \ $ \lambda =  (h_1 )^{p_j}. \prod_{i\not=j} u_i ^{ p_j (\sigma_i +1) - p_i (\sigma_j +1)} $ \ la fonction \ $ \lambda $ \ sera holomorphe au voisinage de l'origine et on aura 
$$  d\lambda \wedge du_j \in u_j .\Omega^{2}_{\mathbb{C}^{2},0} \  . $$
On peut donc \'ecrire \ $ \lambda = \lambda_0  + u_j .L $ \ avec \ $\lambda_0 \in \mathbb{C} $ \ et \ $ L \in \mathcal{O}_{\mathbb{C}^2 , 0 } . $
\begin{itemize}
\item{premier cas }  : il existe \ $ i_o \in [1,k] $ \ tel que l'on ait :\\
 $ p_j (\sigma_{i_0} + 1) - p_{i_0 }(\sigma_j + 1)  >  0  $ \ (on remarquera que le choix de \ $j$ \ fait que tous les entiers \ $ p_j (\sigma_{i_0} + 1) - p_{i_0 }(\sigma_j + 1) $ \ sont positifs ou nuls ).\\
  Alors comme \ $ u_{i_{0}}(0) = 0 $, on aura \ $ \lambda_0 = \lambda (0) = 0 $.  On en d\'eduit que l'on a \ $ h_1 \in  (u_j )$, ce qui prouve le r\'esultat cherch\'e  dans ce cas.
\item{deuxi\`eme cas } :  On a \ $   p_j (\sigma_{i} + 1) - p_{i }(\sigma_j + 1)  =  0, \quad \forall  i \in [1,k]. $\\
On a alors simplement \ $ \lambda = h_1 ^{p_j} $ \ et \ $ d\lambda \wedge du_j \in u_j.\Omega^{2}_{\mathbb{C}^{2},0} $ \  donne 
$$ p_j .h_1 ^{p_j -1}. dh_1 \wedge du_j \in u_j .\Omega^{2}_{\mathbb{C}^{2},0} \  . $$

\noindent Si on a \ $\lambda_0 = 0 $ \ alors l'\'egalit\'e \ $h_1 ^{p_j} = u_j .L $ \ permet imm\'ediatemenet de conclure. Si on a \ $ \lambda_0 \not= 0 $ \ alors \ $ h_1 $ \ est inversible pr\`es de l'origine et on aura
$$  dh_1 \wedge du_j \in u_j .\Omega^{2}_{\mathbb{C}^{2},0} \   $$
d'o\`u  \ $ h_1 = h_1 (0) + u_j .H $ \ avec \ $ H \in \mathcal{O}_{\mathbb{C}^2 ,0 } \ .$

\noindent Posons alors  \ $ \delta := \frac{\sigma_j +1}{p_j} \ ( = \frac{\sigma_i +1}{p_i} \ , \forall i \in [1,k] !) .$ On obtient
\begin{align*}
u_1 ^{\sigma_1} \cdots u_k ^{\sigma_k}.\alpha & = \sum_{l=1}^{k} p_l .u_1 ^{\sigma_1 +1} \cdots u_l ^{\sigma_l} \cdots  u_k ^{\sigma_k +1}.du_l  \\
 \quad         & = \frac{1}{\delta} \sum_{l=1}^{k} (\sigma_l +1).u_1 ^{\sigma_1 +1} \cdots u_l ^{\sigma_l} \cdots  u_k ^{\sigma_k +1}.du_l  =  \frac{1}{\delta} d(u_1 ^{\sigma_1 +1} \cdots u_k ^{\sigma_k +1}).
\end{align*}
On aura donc \ $ d( h_1 (0).u_1 ^{\sigma_1} \cdots u_k ^{\sigma_k}.\alpha ) = 0 \ $ \ et on peut poser \
$ h'_0 : = h_1 (0).u_1 ^{\sigma_1} \cdots u_k ^{\sigma_k}$ \  et \ $ h'_1 = H  $ \ pour avoir \
$  h = h_0 + h'_0 + h'_1 .u_1 ^{\sigma_1} \cdots u_j ^{\sigma_j +1} \cdots  u_k ^{\sigma_k} $ \  avec \\ 
$ d((h_0 + h'_0 ).\alpha) \equiv 0 . \hfill  \blacksquare $
\end{itemize}

\bigskip

\noindent {\bf Remarque.}

\smallskip

\noindent Soit \ $\Delta = p.g.c.d.(p_1 \cdots p_k ). $ Alors si \ $ \Delta = 1 $ \ on ne rencontre pas de \ $ h_0 = u_1 ^{\sigma_1 } \cdots u_k ^{\sigma_k} $ \ v\'erifiant \ $ d(h_0 .\alpha ) = 0$. En effet, si on a \ $ \forall i \in [1,k] \quad \frac{\sigma_i +1}{p_i } = \frac{u}{v} < 1 ,$  avec \ $ p.g.c.d.(u,v) = 1 $, on en d\'eduit que  \ $ \sum_{i=1}^{k} a_i .p_i  = 1 $, donne
\begin{equation*}
\sum_{i=1}^{k}  a_i .p_i .u  =  u  = \sum_{i=1}^{k} a_i .(\sigma_i +1).v
\end{equation*}
ce qui est impossible puisque \ $ u \not\in \mathbb{Z}.v $.

\noindent Par contre, si \ $ \Delta > 1 $, en posant \ $ p_i = \Delta .\sigma_i $ \ on obtient imm\'ediatement que \ $ d( u_1 ^{\sigma_1 -1} \cdots u_k ^{\sigma_k -1}.\alpha) = 0 . $ 

\bigskip

%%%%%%%%%%%%%%%%%%%%%%%%%%%%%%%%%%%%%%%%%%%%%%%%%%%%%%%%%%%%%%%%%%%%%%%%%%%%%%%%%%%%%%%%%%%%%%%%%%%%%%%%
\subsection{ Suspensions.}

\bigskip

\noindent Consid\'erons maintenant une fonction \ $ f : (\mathbb{C}^{n+1},0)  \rightarrow (\mathbb{C},0) $ \ admettant une singularit\'e isol\'ee \`a l'origine et une fonction \ $ g : (\mathbb{C}^{p+1},0)  \rightarrow (\mathbb{C},0) $ \ v\'erifiant l'hypoth\`ese  \ $(HI)$.

\noindent Alors la fonction  \ $ F : (\mathbb{C}^{n+p+1},0)  \rightarrow (\mathbb{C},0) $ \ d\'efinie par \ $ F(x,y) = f(x) + g(y) $ \ v\'erifie l'hypoth\`ese \ $(HI)$  : en effet le lieu singulier de \ $F$ \ est bien de dimension 1 au voisinage de l'origine dans \ $ \mathbb{C}^{n+p+2} $ \ puisqu'il est ensemblistement le produit des lieux singuliers de \ $f$ \ et \ $g$ ;  et si le champ de vecteur holomorphe \ $ \sum_{i=0}^p a_i \frac{\partial}{\partial y_i} $ \ annule \ $g$ \ il annulera \'egalement \ $F$. 

\begin{lemma}\label{KerdF}
\noindent Sous les hypoth\`eses pr\'ecis\'ees ci-dessus on a :
$$  (Ker\, dF)^{n+p+1} = \Omega^{n+1}_{\mathbb{C}^{n+1}} \Join Ker\, dg^p + dF \wedge \Omega^{n+p}_{\mathbb{C}^{n+p+2}} \ .$$
\end{lemma}

\bigskip

\noindent \textit{\underline{Preuve}.} Soient \ $ \pi_1$ \ et \ $ \pi_2 $ \ les projections de  \ $ \mathbb{C}^{n+p+2} \equiv \mathbb{C}^{n+1} \times \mathbb{C}^{p+1} $ \ sur \ $ \mathbb{C}^{n+1}$ \ et \ $ \mathbb{C}^{p+1} $ \ respectivement. Le produit tensoriel externe \ $ \mathcal{F} \Join \mathcal{G} $ \ de faisceaux de \ $\mathcal{O}-$modules (respectivement de complexes de \ $\mathcal{O}-$modules)  sur  \ $ \mathbb{C}^{n+1}$ \ et \ $ \mathbb{C}^{p+1} $ \ d\'esignera le produit tensoriel usuel
$$  \pi_1^*(\mathcal{F}) \otimes   \pi_2^*(\mathcal{G}) \quad {\rm sur} \quad \mathcal{O}_{ \mathbb{C}^{n+p+2}} $$

\noindent Consid\'erons les complexes \ $K^{\bullet}(f) : =  \big( \Omega^{\bullet}_{\mathbb{C}^{n+1}} , \wedge df \big) $ \ et \ $K^{\bullet}(g) : =  \big( \Omega^{\bullet}_{\mathbb{C}^{p+1}} , \wedge dg \big)$ \ o\`u \ $ \Omega^{i}_{\mathbb{C}^{n+1}} $ \ est en degr\'e \ $ i - (n+1) $ \ (resp. \ $ \Omega^{i}_{\mathbb{C}^{p+1}} $ \ en degr\'e \ $ i - (p+1))$.

\noindent Comme \ $f$ \ est \`a singularit\'e isol\'ee en \ $0$ \ dans \ $\mathbb{C}^{n+1} $ \ le complexe \ $ K^{\bullet}(f) $ \ est exact en degr\'es strictement n\'egatifs. Il en est de m\^eme pour \ $ \pi_1^*(K^{\bullet}(f)) $ \ car \ $ \pi_1 $ \ est plat. 

\noindent On a 
$$ K^{\bullet}(F) : = \big(\Omega^{\bullet}_{\mathbb{C}^{n+p+2}} , \wedge dF \big) \simeq K^{\bullet}(f) \Join K^{\bullet}(g)  $$
d'o\`u l'on d\'eduit que 
$$ H^{-1}(K^{\bullet}(F) \simeq H^0 (K^{\bullet}(f)) \Join H^{-1}(K^{\bullet}(g)) \ . $$
Ce qui donne les isomorphismes  :
\begin{align*}
\frac{(Ker\, dF)^{n+p+1}}{dF \wedge \Omega^{n+p}_{\mathbb{C}^{n+p+2}}} \simeq & \frac{ \Omega^{n+1}_{\mathbb{C}^{n+1}}}{df \wedge  \Omega^{n}_{\mathbb{C}^{n+1}}} \Join \frac{(Ker\, dg)^p}{dg \wedge \Omega^{p-1}_{\mathbb{C}^{p+1}}} \\
\smallskip \\
\simeq & \frac{ \Omega^{n+1}_{\mathbb{C}^{n+1}} \Join (Ker\, dg)^p }{ df \wedge  \Omega^{n}_{\mathbb{C}^{n+1}} \wedge \pi_2^*((Ker\, dg)^p) +  \Omega^{n+1}_{\mathbb{C}^{n+1}} \wedge dg \wedge \pi_2^*( \Omega^{p-1}_{\mathbb{C}^{p+1}})} .
\end{align*}
Pour conclure, il suffit de constater que l'on a l'inclusion :
$$df \wedge  \Omega^{n}_{\mathbb{C}^{n+1}} \wedge \pi_2^*((Ker\, dg)^p) +  \Omega^{n+1}_{\mathbb{C}^{n+1}} \wedge dg \wedge \pi_2^*( \Omega^{p-1}_{\mathbb{C}^{p+1}}) \subset dF \wedge \Omega^{n+p}_{\mathbb{C}^{n+p+2}} \quad \quad \blacksquare $$

\begin{prop}\label{T.S.1}

\noindent Sous les hypoth\`eses consid\'er\'ees ci-dessus, si la fonction \ $g$ \ v\'erifie l'hypoth\`ese \ $(HI)$ \ et la condition \ $(P)$, il en sera de m\^eme pour la fonction \ $F$.
\end{prop}

\bigskip

\noindent \textit{Preuve.}  Montrons  que \ $F$ \ v\'erifie la condition \ $(P)$.

\noindent Soit donc \ $ \alpha \in (Ker \, dF)^{n+p+1} $ \ v\'erifiant \ $ d\alpha \in \widehat{J(F)}.\Omega^{n+p+2}_{\mathbb{C}^{n+p+2}} $. Ceci \'equivaut, par d\'efinition de l'id\'eal \ $ \widehat{J(F)}$, \`a demander que \ $d\alpha$ \ soit une section du faisceau \\
 $ dF \wedge \Omega^{n+p+1} $ \ en dehors de l'origine. Mais grace au lemme \ref{KerdF}  on peut \'ecrire
$$  \alpha = dx \wedge \beta + dF \wedge \omega \quad {\rm avec} \quad \beta \in \pi_2^*((Ker \,dg)^p) \ . $$
Comme \ $ d(dF\wedge \omega) \in dF \wedge \Omega^{n+p+1} $ \ la condition impos\'ee \`a \ $\alpha$ \ \'equivaut \`a demander que \ $ dx \wedge \beta $ \ soit dans \ $ \big(J(f) + J(g)\big).\Omega^{n+p+2} $ \ en dehors de l'origine. D\'ecomposons alors
$$ \mathcal{O}_{\mathbb{C}^{n+1},0} = J(f)_0 \oplus V $$
o\`u \ $V$ \ est un \ $\mathbb{C}-$espace vectoriel de dimension finie (\'egale au nombre de Milnor de \ $f$ \ en \ $0$). D\'ecomposons le germe de \ $ \beta$ \ le long de \ $ \lbrace 0 \rbrace \times \mathbb{C}^{p+1} $ \ sous la forme \ $\beta = \beta_0 + \beta_1 \quad $ \ avec 
$$  \beta_0 \in J(f)_0 \Join (Ker \,dg)^p \quad {\rm et} \quad  \beta_1 \in V \otimes_{\mathbb{C}} \pi_2^{-1}((Ker \,dg)^p) $$ 
o\`u  \ $\pi_2^{-1} $ \ d\'esigne l'image r\'eciproque ensembliste du faisceau .

\noindent Comme la diff\'erentielle partielle en \ $y$ \ laisse \ $ J(f) $ \ stable, la condition sur \ $d\alpha$ \ est \'equivalente \`a demander que \ $ d_{/y}\beta_1 \in V \otimes J(g).\Omega^{p+1}_{\mathbb{C}^{p+1}} $ \ en dehors de l'origine dans \ $ \mathbb{C}^{p+1} $ \footnote{ on remarquera que \ $ J(f).dx \Join Ker \,dg^p \subset dF\wedge \Omega^{n+p+1}$ \ et donc que l'on a  
$$ \alpha = dx \wedge \beta_1 \quad  {\rm modulo} \quad dF\wedge \Omega^{n+p} \ .$$ }.
Mais puisque \ $g$ \ v\'erifie la propri\'et\'e \ $(P)$ \ on peut \'ecrire
\begin{align*}
 \beta_1 = \Lambda + dg \wedge M  & \quad {\rm avec} \quad \Lambda \in V \otimes \pi_2^{-1}\big[(Ker\, dg_0)^p \cap Ker\, d\big] \\
\smallskip \\
\quad & {\rm et} \quad M \in V \otimes \pi_2^{-1}\big[\Omega^{p-1}_{\mathbb{C}^{p+1},0} \big]
\end{align*}
On conclut alors car on a 
\begin{align*}
dx \wedge \Lambda\  \in \ \Omega^{n+1}_{\mathbb{C}^{n+1},0} \Join \big((Ker\, dg)^p \cap Ker\, d \big)_0 \subset \big((Ker\, dF)^{n+p+1} \cap Ker \, d \big)_0 \\
{\rm et} \qquad dx \wedge dg \wedge M \in \big(dF \wedge  \Omega^{n+p}\big)_0 \qquad \qquad \qquad \blacksquare
\end{align*}

\bigskip

\noindent  {\bf Remarque.}

\bigskip

\noindent La proposition pr\'ec\'edente combin\'ee avec les propositions \ref{a.b.s.} et \ref{s.b.s.} fournit beaucoup d'exemples de fonctions v\'erifiant l'hypoth\`ese \ $(HI)$ \ et satisfaisant \'egalement la condition  \ $(P)$ . 

\noindent Plus pr\'ecis\'ement, elle donne, pour chaque nouvel exemple de  fonction v\'erifiant  \ $(HI)$ \ et satisfaisant la condition  \ $(P)$ \ une nouvelle famille d'exemples. Cela justifie d'explorer la condition \ $(P)$ \ sur les fonctions les plus simples v\'erifiant \ $(HI)$ \ pour \ $n=2$.

\begin{prop}\label{T.S.2}
\noindent On suppose que la fonction\ $f$ \ est  \`a singularit\'e isol\'ee \`a l'origine de \ $\mathbb{C}^{n+1}$ \  et que la fonction \ $g$ \ v\'erifie l'hypoth\`ese \ $(HI)$ \ et la condition \ $(P)$, au voisinage de l'origine dans \ $\mathbb{C}^{p+1}$. Consid\'erons les (pre)-(a,b)-modules
\begin{equation*}
E'_f : = \frac{ \Omega^{n+1}_{\mathbb{C}^{n+1},0}}{df \wedge d\Omega^{n-1}_{\mathbb{C}^{n+1},0}} \ \ 
E'_g : = \frac{ \Omega^{p+1}_{\mathbb{C}^{p+1},0}}{d((Ker\, dg)^p_0)}
 \quad {\rm et}\quad  E'_F : = \frac{ \Omega^{n +p+2}_{\mathbb{C}^{n+p+2},0}}{d((Ker\, dF)^{n+p+1}_0)} 
\end{equation*}

\noindent Alors l'application donn\'ee  par produit ext\'erieur
$$ \Lambda : E'_f {\otimes}_{a,b} E'_g  \to E'_F$$

induit un isomorphisme entre les (a,b)-modules associ\'es.
\end{prop}

\bigskip

\noindent \textit {\underline{Preuve}.} Gr\^ace au lemme de Nakayama, il suffit de voir que l'application induite \ $ E'_f/b.E'_f \otimes E'_g/b.E'_g \to E'_F/b.E'_F $ \ est un isomorphisme. Ceci r\'esulte facilement des isomorphismes suivants
$$ \mathcal{O}/J(f) \otimes \widehat{J(g)}/J(g) \simeq \frac{\widehat{J(g)}}{J(f)\cap\widehat{J(g)} + J(g)}\simeq \widehat{J(F)}/J(F) $$
$$ \mathcal{O}/J(f)\otimes E'_g/jE''_g \simeq  E'_F/jE''_F $$
et du diagramme commutatif aux lignes exactes :\\

\smallskip

\xymatrix{ 0\ar[r] &  \mathcal{O}/J(f) \otimes \widehat{J(g)}/J(g)\ar[d] \ar[r] & E'_f/b.E'_f \otimes E'_g/b.E'_g \ar[d] \ar[r] & \mathcal{O}/J(f)\otimes E'_g/jE''_g \ar[d] \ar[r] & 0 \\
0 \ar[r] & \widehat{J(F)}/J(F) \ar[r] & E'_F/b.E'_F \ar[r] &  E'_F/jE''_F \ar[r] & 0 .}

\bigskip

Ceci ach\`eve la d\'emonstration. $ \hfill \blacksquare$

\bigskip

Pour une version un peu plus sophistiqu\'ee de ce r\'esultat \`a la "Thom-Sebastiani", \\
le lecteur pourra consulter [B.S. 04].

\section*{References}

\bigskip
\begin{itemize}
\item{[B.84]} D.Barlet : {\it Contribution du cup-produit de la fibre de Milnor aux poles de  ...}
Ann. Inst. Fourier (1984) t. 34, fasc.4, p.75-106.
\item{[B.91]} D. Barlet : {\it Interaction de strates cons{\'e}cutives pour les cycles
             {\'e}vanescents.} Ann. Scient. Ec. Norm. Sup., 4{\'e}me s{\'e}rie, t. 24, (1991), p. 401-506.
\item{[B.93]} D. Barlet : {\it Th{\'e}orie des $\;(a,b)-$modules I.} in Complex Analysis
             and Geometry, Plenum Press, (1993), p.1-43.
\item{[B.95]} D. Barlet : {\it Th{\'e}orie des $\;(a,b)-$modules II. Extensions.} in
             Complex Analysis and Geometry, Pitman Research, Notes in Mathematics, s{\'e}ries 366,
             (1997),  p.19-59.
\item{[B.05]} D. Barlet : {\it Interaction de strates cons\'ecutives II .}  Publ. RIMS Kyoto University, vol 41 (2005), p.139-173.
\item{[B.04 a]} D. Barlet : {\it Interactions de strates cons\'ecutives ... III : le cas de la valeur propre 1.}  Pr\'epublication de l' Institut E. Cartan (Nancy), 2004/38
\item{[B.04 b]} D. Barlet : {\it Sur certaines singularit\'es non isol\'ees d'hypersurfaces II.} Pr\'epublication de l' Institut E. Cartan (Nancy), 2005/42
\item{[Be.01]} R. Belgrade : {\it Dualit\'e et Spectres des (a,b)-modules.} Journal of Algebra 245, (2001), p.193-224.
\item{[Bj.93]} J.-E. Bjork : {\it Analytic D-modules and Applications.} Mathematics and its Applications, vol.247, Kluwer 1993. 
\item{[Br.70]} E. Brieskorn : {\it Die Monodromie der isolierten Singularit\"aten von Hyperfl\"achen.} Manuscripta Math. 2 (1970), p.103-161.
\item{[B.S. 04]} D. Barlet and M. Saito : {\it Brieskorn Modules and Gauss-Manin systems for non isolated hypersurfaces singularities.} Pr\'epublication de l' Institut E. Cartan (Nancy), 2004/54.
\item{[G.66]} A. Grothendieck : {\it On the de Rham cohomology of algebraic varieties.} Publ. Math. IHES 29 (1966), p.93-101.
\item{[H.64]} H. Hironaka : {\it Resolution of singularities of an algebraic variety over a field of characteristic zero I. II. } Ann. Math. 79 (1964), p.109-203  and p. 205-326.
\item{[K.75]} M. Kashiwara : {\it On the maximally over determined systems of differential
             equations.} Publ. R.I.M.S., vol.10, (1975), p.563-579.
\item{[M.74]} B.Malgrange : {\it Int\'egrale asymptotique et monodromie.} Ann. Scient. Ec. Norm. Sup. , t.7 (1974), p.405-430 ( On pourra consulter l'Appendice de [B.84] pour des d\'etails sur le th\'eor\`eme de positivit\'e pour un germe de fonction holomorphe (r\'eduite)  arbitraire.)
\item{[S.70]} M. S\'ebastiani : {\it Preuve d'une conjecture de Brieskorn.} Manuscripta Math. 2 (1970), p.301-308.

\end{itemize}

\bigskip

\noindent Daniel Barlet, 

\noindent Universit\'e Henri Poincar\'e (Nancy I ) et Institut Universitaire de France,

\noindent Institut E.Cartan  UHP/CNRS/INRIA, UMR 7502 ,

\noindent Facult\'e des Sciences et Techniques, B.P. 239

\noindent 54506 Vandoeuvre-les-Nancy Cedex , France.

\noindent e-mail :  barlet@iecn.u-nancy.fr

\end{document}